\documentclass[10pt]{article}
\usepackage{amssymb,amsmath,amsthm}
\usepackage{graphicx}
\usepackage{epstopdf}
\usepackage{exscale}
\usepackage{relsize}
\usepackage{subfigure}
\newcommand{\R}{\mathbb{R}}

\begin{document}
\title{\Large\bf{Existence of solutions for a quasilinear elliptic system with local nonlinearity on $\R^N$}
 }
\date{}
\author {  Xingyong Zhang$^{1,2}$\footnote{Corresponding author, E-mail address: zhangxingyong1@163.com}, Cuiling Liu$^{1}$  \\
      {\footnotesize $^1$Faculty of Science, Kunming University of Science and Technology, Kunming, Yunnan, 650500, P.R. China.}\\
        {\footnotesize  $^2$School of Mathematics and Statistics, Central South University,
        Changsha, Hunan, 410083, P.R. China.}\\}
 \date{}
 \maketitle

 \begin{center}
 \begin{minipage}{15cm}
 \par
 \small  {\bf Abstract:} In this paper, we investigate the existence of solutions for a class of quasilinear elliptic system
   \begin{eqnarray*}\label{aa1}
 \begin{cases}
  -\mbox{div}(\phi_1(|\nabla u|)\nabla u)+V_1(x)\phi_1(|u|)u=\lambda F_u(x, u,v), \ \ x\in \R^N,\\
   -\mbox{div}(\phi_2(|\nabla v|)\nabla v)+V_2(x)\phi_2(|v|)v=\lambda F_v(x, u,v), \ \ x\in \R^N,\\
  u\in W^{1,\Phi_1}(\R^N), v\in W^{1,\Phi_2}(\R^N),
   \end{cases}
 \end{eqnarray*}
where $N\ge 2$, $\inf_{\R^N}V_i(x)>0,i=1,2$, and $\lambda>0$. We obtain that when the nonlinear term $F$ satisfies some growth conditions only in a circle with center $0$ and radius $4$, system has a nontrivial solution $(u_\lambda,v_\lambda)$ with  $\|(u_{\lambda},v_{\lambda})\|_{\infty}\le 2$  for every $\lambda$ large enough, and the families of solutions $\{(u_\lambda,v_\lambda)\}$ satisfy that  $\|(u_\lambda,v_\lambda)\|\to 0$ as $\lambda\to \infty$. Moreover, a corresponding result for a quasilinear elliptic equation is also obtained, which is better than the result for the elliptic system.

 \end{minipage}
 \end{center}
  \allowdisplaybreaks
 \vskip2mm
 {\section{Introduction}}
 \setcounter{equation}{0}
 \par
 In the past years, the existence and multiplicity of solutions for the quasilinear elliptic problem with the following form
 \begin{eqnarray}\label{bb1}
\begin{cases}
  -\mbox{div}(\phi(|\nabla u|)\nabla u)+V(x)\phi(|u|)u=f(x,u), \ \ x\in \Omega,\\
  u=0,x\in \partial \Omega
   \end{cases}
 \end{eqnarray}
 has been investigated extensively (for example,
 see \cite{Alves2014},  \cite{Mahiout2020}, \cite{Bonanno201201}, \cite{Bonanno201202}, \cite{Fukagai1995},  \cite{Fukagai2006},
 \cite{Fukagai2009}, \cite{Liu2019},
 \cite{Mihailescu2011} and references therein), where   $\Omega\subset \R^N$ is an open set, $N\ge 2$, $V,f$ are continuous functions and $\phi:[0,\infty)\to [0,\infty)$ satisfies some suitable monotonicity and growth conditions.
Equation (\ref{bb1}) arises from some fields of physics, for example,\\
  $(1)$ nonlinear elasticity: $\Phi(t)=(1+t^2)^\gamma-1, \gamma>\frac{1}{2};$\\
 $(2)$ plasticity: $\Phi(t)=t^\alpha(\log(1+t))^\beta, \alpha\geq 1, \beta>0;$\\
 $(3)$ generalized Newtonian fluids: $\Phi(t)=\int_{0}^{t}s^{1-\alpha}(\sinh^{-1}s)^\beta ds, 0\leq \alpha \leq 1, \beta>0,$\\
where $\Phi(t)=\int_0^t \phi(s)sds$ (see \cite{Alves2014}, \cite{Fuchs199801}, \cite{Fuchs02}, \cite{Fukagai1995}, \cite{Fukagai2006}).

\par
 Specifically,
in \cite{Alves2014},  Alves-Figueiredo-Santos considered equation (\ref{bb1}) with $\Omega= \R^N$.
They assumed that $V\in C(\R^N,\R)$, $\inf_{x\in \R^N} V(x)=V_0>0$, $V$ is a radial function or is a $\mathbb Z^N$ periodic function and $f\in C(\R^N,\R)$ and satisfies
 \begin{eqnarray}\label{aaa2}
  \lim_{|t|\to 0} \frac{f(t)}{\phi(|t|)|t|}=0,\quad \lim_{|t|\to +\infty}\frac{f(t)}{\phi_*(|t|)|t|}=0
 \end{eqnarray}
 and there exists $\nu>m$ such that
 \begin{eqnarray}\label{aaa3}
 0<\nu F(t)=\int_0^t f(s)ds\le tf(t)\ \ \mbox{ for all } t\in \R/\{0\}
\end{eqnarray}
 which is usually called as Ambrosseti-Rabinowitz condition ((AR) for short). After developing a Strauss-type result and a Lions-type result,  they obtained equation (\ref{bb1}) has a nontrivial solution.
Recently,
in \cite{Liu2019}, Liu considered the case that $V$ has an infinite potential well, that is,
\par
$(\mathcal{V}1)$ for all $M>0$, $\mu(V^{-1}(-\infty, M])<\infty$, where $\mu$ is the Lebesgue measure,\\
  or $V$ has  a finite potential well, that is,
  \par
   $(\mathcal{V}1)'$ for all $x\in \R^N$, $V(x)<\lim_{|x|\to \infty}V(x)<\infty$.\\
   He also considered the case that $V$ is a steep potential well, that is, $V(x)=\lambda a(x)+1$, where $\lambda$ is  a parameter and $a\in C(\R^N,\R)$. For all these cases, he assumed that (\ref{aaa2}) and (AR) hold. Then he obtained some existence and multiplicity results of solutions for system (\ref{bb1}).
\par
It is easy to see that (\ref{aaa2}) and (AR) imply that $f$ satisfies some conditions near  both  $0$ and  $\infty$.  So it is natural to ask if it is possible to restrict those conditions for $f$ to either of them.
To this end,
in \cite{costa},  Costa and Wang  used a cut-off technique together with energy estimates to study the multiplicity of both signed and sign-changing solutions for one-parameter family of elliptic problems (\ref{bb1}) with $\phi=1,\; V=0,\; f(x,u)=\lambda f(u),$ where
 $\lambda>0$ is a parameter, $\Omega$ is a bounded smooth domain in $\R^N (N \ge 3)$ and $f\in C^1(\R,\R)$. By such  ingenious method, the nonlinearity $f(u)$ was assumed to satisfy the superlinear growth only in a neighborhood of $u = 0$.
Afterwards, in \cite{Medeiros},   Medeiros and Severo applied the idea in \cite{costa} to the problem (\ref{bb1}) with $\phi(t)=|t|^{p-2}$ and $  f(x,u)=\lambda f(u)$ on the whole space $\R^N$, i.e. the following $p$-Laplacian equation
  \begin{equation}\label{kkkkk3}
-\Delta_p u +V(x)|u|^{p-2}u= \lambda f(u) \ \ \mbox{ in } \R^N,
 \end{equation}
where $1<p<N$ and $\lambda>0$. They assumed that $V\in C(\R^N,\R)$, $\inf_{x\in \R^N} V(x)=V_0>0$ and (V1) holds. Moreover, $f$ satisfies the following conditions:\\
$(f1)$ there exists $r\in (p,p^*)$ such that
$$
\limsup_{|s|\to 0}\frac{f(s)s}{|s|^r}<+\infty;
$$
$(f2)$ there exists $q\in (p,p^*)$ such that
$$
 \liminf_{|s|\to 0}\frac{F(s)}{|s|^q}>0;
$$
$(f3)$ there exists $\nu\in (p,p^*)$ such that
$$
0<\nu F(s)\le sf(s)\mbox{ for } |s|\not =0 \mbox{ small} ,
$$
where $p^*=\frac{Np}{N-p}$. With developing Moser iteration technique,  they proved that equation (\ref{kkkkk3}) has one positive solution and one negative solution for all $\lambda$ large enough. $(f1)$-$(f3)$ show that  $f(s)$ satisfies the superlinear growth only in a neighborhood of $s = 0$.   Here, it needs to be emphasized that $(f1)$-$(f3)$ with $p=2$ were given first in \cite{costa}. The idea in \cite{costa} has been applied to various differential equations and we cite \cite{Medeiros2008}, \cite{Huang2019} and \cite{Xu2016}  as some examples.

 \par
Inspired by \cite{costa} and \cite{Medeiros}, in this paper, we investigate the existence of solutions for the following quasilinear elliptic problem with a parameter
   \begin{eqnarray}\label{aa1}
\begin{cases}
  -\mbox{div}(\phi_1(|\nabla u|)\nabla u)+V_1(x)\phi_1(|u|)u=\lambda F_u(x, u,v), \ \ x\in \R^N,\\
   -\mbox{div}(\phi_2(|\nabla v|)\nabla v)+V_2(x)\phi_2(|v|)v=\lambda F_v(x, u,v), \ \ x\in \R^N,\\
  u\in W^{1,\Phi_1}(\R^N), v\in W^{1,\Phi_2}(\R^N),
   \end{cases}
 \end{eqnarray}
where $N\ge 2$, $\lambda$ is a parameter with $\lambda>0$, $F\in C^1(\R^N\times \R\times\R,\R)$,    $\phi_i$ and $V_i\in C(\R^N,\R^+)$, $i=1,2$ satisfy the following assumptions:\\
$(\phi_1)$ $\phi_i\in C^1(0,\infty)$,  $t\rightarrow \phi_i(t)t$ are strictly increasing, $i=1,2$;\\
 $(\phi_2)$ $1<l_i:=\inf_{t>0}\frac{t^2\phi_i(t)}{\Phi_i(t)}\leq\sup_{t>0}\frac{t^2\phi_i(t)}{\Phi_i(t)}=:m_i<\min\{N, l_i^*\}$, where $\Phi_i(t):=\int_{0}^{|t|}s\phi_i(s)ds, t\in \mathbb{R}$ and $l_i^*:=\frac{l_iN}{N-l}$, $i=1,2$.
  \par
Recently, in \cite{Wang2016} and \cite{Wang2017}, Wang-Zhang-Fang investigated the
quasilinear elliptic system (\ref{aa1}) with $\lambda=1$.
 In \cite{Wang2016}, when $F$ has a sub-linear growth, by using the least
action principle, they obtained that system has at least one nontrivial solution and when $F$ satisfies
an additional symmetric condition, by using the genus theory, they obtained that system has
infinitely many solutions. In \cite{Wang2017}, by using the mountain pass theorem, when $F$ satisfies some superlinear conditions, they obtained that system has a ground state solution. Especially, they obtained the following theorem.
\vskip2mm
\noindent
{\bf Theorem A.} (\cite{Wang2017}) {\it Assume that $(\phi_1)$, $(\phi_2)$ and the following conditions hold:\\
(V1)\  $V_i$ are $1$-periodic functions in $x_1,\cdots,x_N$  for all $x\in \R^N$, $i=1,2$ (called $1$-periodic for short);\\
(V2)\ there exit two positive constants $\alpha_1$ and $\alpha_2$ such that
$$
\alpha_1\le \min\{V_1(x),V_2(x)\}\le \max\{V_1(x),V_2(x)\}\le \alpha_2
$$
for all $x\in \R^N$;\\
(H0) $F\in C^1(\R^N\times \R\times \R,\R)$, $F$ is $1-$periodic in $x\in \R^N$ and $F(x,0,0)=0$ for all $x\in \R^N$;\\
(H1)
\begin{eqnarray*}
&  &  \lim_{|(t,s)|\to 0}\frac{F_t(x,t,s)}{\phi_1(|t|)+{\widetilde{\Phi}_1}^{-1}(\Phi_2(|s|))}=0,\ \ \lim_{|(t,s)|\to 0}\frac{F_s(x,t,s)}{{\widetilde{\Phi}_2}^{-1}(\Phi_1(|t|))+\phi_2(|s|)}=0,\\
&  & \lim_{|(t,s)|\to \infty}\frac{F_t(x,t,s)}{\Phi_{1*}'(|t|)+{\widetilde{\Phi}_{1*}}^{-1}(\Phi_{2*}(|s|))}=0,\ \ \lim_{|(t,s)|\to \infty}\frac{F_s(x,t,s)}{{\widetilde{\Phi}_{2*}}^{-1}(\Phi_{1*}(|t|))+\Phi_{2*}'(|s|)}=0\\
\end{eqnarray*}
uniformly in $x\in \R^N$, where $\widetilde{\Phi}_i$ and $\Phi_{i*}$, $i=1,2$ are defined in section 2 below;\\
(H2) there exist $\mu_i>m_i(i=1,2)$ such that
$$
0<F(x,t,s)\le \frac{1}{\mu_1}tF_t(x,t,s)+\frac{1}{\mu_2}sF_s(x,t,s),\ \ \mbox{for all } (t,s)\not=0.
$$
Then system (\ref{aa1}) with $\lambda=1$ has a ground state solution.}
\par
 It is easy to see (H0)-(H2) imply  $F$ satisfies some growth near both $0$ and $\infty$. In this paper, by applying the method in \cite{costa} and developing the Moser iteration technique, we only need to make some assumptions on the nonlinearity $F(x,t,s)$ in a circle with center $0$ and radius $4$. If we assume that $V_i,i=1,2$ satisfies $(\mathcal{V}1)$ instead of $({V}1)$ in Theorem 1.1 below, then Theorem 1.1 can be seen as an extension of the result in \cite{Medeiros} to system (\ref{aa1}) in some sense. Moreover, we shall find that the elliptic system  (\ref{aa1})  is more complex and more general than the scalar equation  (\ref{bb1}) with $x\in \R^N$, which directly leads to some  stronger restrictions for the nonlinearity $F$ in Theorem 1.1 because of a different Moser iteration result, see section 5 for more details.  To be precise, we obtain the following theorem for system (\ref{aa1}):
\vskip2mm
\noindent
{\bf Theorem 1.1.} {\it Assume that $(\phi_1)$-$(\phi_2)$, (V1) and the following conditions hold:\\
 $(\phi_3)$ there exist positive constants $q_i$ such that $t^2\phi_i(|t|)\ge q_i|t|^{l_i}$ for all  $t\in \R$, $i=1,2$;\\
 (V0)\ $V_i\in C(\R^N,\R^+)$ and $V_{i,\infty}:=\inf_{\R^N}V_i(x)>0$, $i=1,2$;\\
 (F0) $F\in C^1(\R^N\times \R\times\R,\R)$ and $F$ is $1-$periodic in $x\in \R^N$;\\
(F1)\ there exist $k_i\in (m_i,l_i^*)$ and  $M_i>0$, $i=1,2$  such that
$$
F(x, t,s)\ge M_1 |t|^{k_1}+M_2 |s|^{k_2}
$$
for all $|(t,s)|<4$  and $x\in  \R^N$;\\
(F2)\ there exist  $M_3>0,M_4>0$, $r_1\in \bigg(\max\left\{m_1,m_2,1+\frac{(1+m_1)(l_2-1)(\Theta_2-1)}{l_2\Theta_2}\right\},\min\left\{l_1^*,
1+ \frac{l_1^*(l_1-1)}{l_1}-\frac{(l_1-1)(1+m_1)}{\Theta_1 l_1}\right\}\bigg)$ and
$r_2\in \bigg(\max\left\{m_1,m_2,1+\frac{(1+m_2)(l_1-1)(\Theta_1-1)}{l_1\Theta_1}\right\},\min\left\{l_2^*,
1+ \frac{l_2^*(l_2-1)}{l_2}-\frac{(l_2-1)(1+m_2)}{\Theta_2 l_2}\right\}\bigg)$  for some  $\Theta_1>1$ and $\Theta_2>1$,  such that
\begin{eqnarray*}
   |F_t(x,t,s)|\le M_3 |t|^{r_1-1}+M_4 |s|^{r_2-1},\ \  |F_s(x,t,s)|\le M_3 |t|^{r_1-1}+M_4 |s|^{r_2-1}
\end{eqnarray*}
for all  $|(t,s)|< 4$  and $x\in  \R^N$;\\
(F3) there exist $\mu_i>m_i$,  $i=1,2$ such that
$$
0< F(x,t,s)\le \frac{1}{\mu_1}tF_t(x,t,s)+\frac{1}{\mu_2}sF_s(x,t,s)
$$
for all $x\in  \R^N$ and  $|(t,s)|<4$  with $(t,s)\not=(0,0)$.\\
Then there exists $\Lambda_*>0$ such that system (\ref{aa1})  has a nontrivial solution $(u_\lambda,v_\lambda)$ with   $\|(u_{\lambda},v_{\lambda})\|_{\infty}\le 2$ for each $\lambda> \Lambda_*$ and $\|(u_\lambda,v_\lambda)\|\to 0$ as $\lambda\to\infty$.}

\vskip2mm
\noindent
{\bf Remark 1.1.}  There exist examples satisfying Theorem 1.1. For example, let
$$F(x,t,s)= \sigma(t,s)b(x)G_1(t,s)+(1-\sigma(t,s))b(x)G_2(t,s)$$
 for all $(x,t,s)\in \R^{N}\times \R\times \R$, where $b(x)$ satisfies (V0), $G_1\in C^1(\overline{B_8},\R)$ satisfies (F1)-(F3), $G_2\in C^1(\R^2/B_4,\R)$, $B_R$ denotes a circle with center $0$ and radius $R$,  and
   \begin{eqnarray}\label{dd2}
 \sigma(t,s)= \begin{cases}
           1, \;\;\; \;\;\; \;\;\; \;\;\; \;\;\; \;\;\; \;\;\; \;\;\; \;\;\; \;\;\; \;\;\; \; \text { if }\;\;|(t,s)| < 4, \\
             \sin\dfrac{\pi(t^2+s^2-64)^2}{4608}, \;\;\;  \text { if }\;\;4 \le |(t,s)| \leq 8, \\
           0, \;\;\; \;\;\; \;\;\; \;\;\; \;\;\; \;\;\; \;\;\; \;\;\; \;\;\; \;\;\; \;\;\; \;  \text { if }\;\;|(t,s)| > 8.
          \end{cases}
   \end{eqnarray}
  In particular, $F$ satisfies Theorem 1.1 but not satisfying Theorem A if we let $N=6$, $\phi_1(t)=4|t|^{2}+5|t|^{3}$ and $\phi_2(t)=4|t|^{2}\log(2+| t|)+\dfrac{|t|^{3}}{1+|t|}$, $b(x)= \left(1+\sum_{i=1}^{6}\cos^2\pi x_i\right)$ (or $b(x)\equiv 1$), $G_1(t,x)=|t|^{\frac{17}{2}}+|s|^{\frac{17}{2}}+|t|^{7}|s|^{7}$ and $G_2(t,s)=|t|^3+|s|^3$.
The readers can see  the detail computation  in section 4.
\vskip2mm
\noindent
{\bf Remark 1.2.}  In Theorem 1.1, we only assume that $V_i(i=1,2)$ are periodic functions. In fact, it is also possible that similar results can be established if  $V_i(i=1,2)$  are radial functions or satisfy $(\mathcal{V}1)$ or $(\mathcal{V}1)'$ by combing those arguments in \cite{Alves2014} and \cite{Liu2019}.
\vskip2mm
\par
We organize the paper as follows. In section 2, we recall some knowledge for Orlicz and Orlicz-Sobolev spaces.  In section 3, we complete the proof of Theorem 1.1. In section 4, we present some detail arguments for the example mentioned in Remark 1.1. In section 5, corresponding to Theorem 1.1, we present a result for a quasilinear elliptic equation with a parameter $\lambda$, which shows some differences  between equation (\ref{bb1}) with $\Omega=\R^N$ and system (\ref{aa1}).
 \vskip2mm
 {\section{Preliminaries}}
  \setcounter{equation}{0}
 \par
In this section, we recall some notions and  properties about Orlicz and Orlicz-Sobolev spaces and some useful lemmas.
  The readers can see these details in \cite{Adams2003, Alves2014, Brezis1991, Fukagai2006, M.A. Krasnoselski1961, Rabinowitz1986, M. M. Rao2002, Willem M1997}.
  \par
We will start with some properties  about $N$-function. Assume that $a: [0,\infty)\rightarrow [0,\infty)$  is a right continuous, monotone increasing function with \\
 (i) $a(0)=0$;\\
 (ii) $\lim_{t\rightarrow \infty} a(t)=\infty$;\\
 (iii) $a(t)>0 $ whenever $t>0$.\\
 Then the integral $A(t)=\int_{0}^{t}a(s)ds$
 is called an $N$-function which is defined on $[0,\infty)$.
  \par
    Define the complement of $A$  by
 $$\widetilde{A}(t)=\max_{s\geq 0}\{ts-A(s)\},\quad\mbox{ for } t\geq 0.$$
 Then $\widetilde{A}$ is also an $N$-function and $\widetilde{\widetilde{A}}=A$,
 and the Young's inequality holds, that is
 \begin{eqnarray}\label{u1}
 st\leq A (s)+\widetilde{A}(t), \quad\mbox{ for all } s, t\geq 0,
 \end{eqnarray}
  If
 $$
 \sup_{t>0}\frac{A(2t)}{A(t)}<\infty,
 $$
 then we call $A$ satisfies a $\Delta_2$-condition globally.
 When $A$ satisfies $\Delta_2$-condition globally,
 \begin{eqnarray}\label{2.2tj}
A(t)>A^{\beta}(t)
\end{eqnarray}
 for any $\beta > 1$ (see \cite{M.A. Krasnoselski1961}) and the Orlicz space $L^{A}(\Omega)$ is defined by  the vectorial space consisting of the measurable functions $u: \Omega\rightarrow \mathbb{R}$ satisfying
$$\int_{\Omega}A(|u|)dx<\infty,$$
where $\Omega \subset \mathbb{R}^N$ is an open set. Define
$$\|u\|_{A}:=\inf \left\{\alpha>0: \int_{\Omega}A \left(\frac{|u|}{\alpha}\right)dx<1\right\},\quad \mbox{ for } u\in L^{A}(\Omega), $$
which is called  Luxemburg norm. Then  $(L^{A}(\Omega),\|\cdot\|_A)$ is a Banach space.
If $A(t)=|t|^p (1<p<+\infty)$,   $(L^{A}(\Omega),\|\cdot\|_A)$  corresponds to the classical Lebesgue space $L^p(\Omega)$ with the norm
$$\|u\|_{p}:=\left(\int_{\Omega}|u(x)|^pdx\right)^{\frac{1}{p}}.$$

\par
Define
$$W^{1,A}(\Omega)=\left\{u \in L^{A}(\Omega): \frac{\partial u}{\partial x_i} \in L^{A}(\Omega), i=1,\cdots, N\right\}$$
 with the norm
$$\|u\|_{1, A}=\|u\|_{A}+\|\nabla u \|_{A}.$$
Then $W^{1, A}(\Omega)$ is a Banach space which is called Orlicz-Sobolev space. Denote the closure of $C_{0}^{\infty}(\Omega)$ in $W^{1, A}(\Omega)$ by $W_{0}^{1, A}(\Omega)$. $W_{0}^{1, A}(\mathbb{R}^N)=W^{1, A}(\mathbb{R}^N)$ if $\Omega=\mathbb{R}^{N}$.
\vskip2mm
 \noindent
 {\bf Lemma 2.1.} (\cite{Fukagai2006})\ \  {\it If $A$ is an $N$-function, then the following conditions are equivalent:\\
(i)
\begin{eqnarray}\label{ 2.4+}
l=\inf_{t>0}\frac{ta(t)}{A(t)}\leq\sup_{t>0}\frac{ta(t)}{A(t)}=m;
\end{eqnarray}
(ii) let $\zeta_0(t)=\min\{t^l, t^m\}$, $\zeta_1(t)=\max\{t^l, t^m\}$, for $t\geq0$. $A$ satisfies
$$\zeta_0(t)A(\rho)\leq A(\rho t)\leq \zeta_1(t)A(\rho), \quad \mbox{ for all } \rho, t\geq 0;$$
(iii) $A$ satisfies a $\Delta_2$-condition globally.
}
\vskip2mm
\noindent
 {\bf Lemma 2.2.} (\cite{Fukagai2006})\ \  {\it If $A$ is an $N$-function and (\ref{ 2.4+}) holds, then $A$ satisfies }
$$\zeta_0(\|u\|_{A})\leq\int_{\mathbb{R}^N}A(|u|)dx\leq\zeta_1(\|u\|_{A}), \quad \mbox{ \it for all }u \in L^{A}(\mathbb{R}^N).$$
 \vskip2mm
\noindent
 {\bf Lemma 2.3.} (\cite{Fukagai2006})\ \  {\it If $A$ is an $N$-function, $l$, $m\in(1, \infty )$ and (\ref{ 2.4+}) holds. Let $\widetilde{A}$ be the complement of $A$ and $\zeta_2(t)=\min\{t^{\widetilde{l}},t^{\widetilde{m}}\}$, $\zeta_3(t)=\max\{t^{\widetilde{l}},t^{\widetilde{m}}\}$, for $t\geq0$, where $\widetilde{l}:=\frac{l}{l-1}$ and $\widetilde{m}:=\frac{m}{m-1}$. Then $\widetilde{A}$ satisfies\\}
{\it (i)
$$\widetilde{m}=\inf_{t>0}\frac{t\widetilde{A}^{\prime}(t)}{\widetilde{A}(t)}\leq\sup_{t>0}\frac{t\widetilde{A}^{\prime}(t)}{\widetilde{A}(t)}=\widetilde{l};$$
(ii)
$$\zeta_2(t)\widetilde{A}(\rho)\leq \widetilde{A}(\rho t)\leq \zeta_3(t)\widetilde{A}(\rho), \quad \mbox{\it for all } \rho, t\geq 0;$$
(iii)}
$$\zeta_2(\|u\|_{\widetilde{A}})\leq\int_{\mathbb{R}^N}\widetilde{A}(|u|)dx\leq\zeta_3(\|u\|_{\widetilde{A}}), \quad \mbox{\it for all }u \in L^{\widetilde{A}}(\mathbb{R}^N).$$
\vskip2mm
\noindent
 {\bf Lemma 2.4.} (\cite{Fukagai2006})\ \  {\it If $A$ is an $N$-function, $l$, $m\in(1, N)$ and (\ref{ 2.4+}) holds. Let $\zeta_4(t)=\min\{t^{l^*},t^{m^*}\}$, $\zeta_5(t)=\max\{t^{l^*},t^{m^*}\}$, for $t\geq0$, where $l^*:=\frac{lN}{N-l}$, $m^*:=\frac{mN}{N-m}$. Then $A_{*}$ satisfies\\
(i)
$$l^*=\inf_{t>0}\frac{tA_{*}^{\prime}(t)}{A_*(t)}\leq\sup_{t>0}\frac{tA_{*}^{\prime}(t)}{A_*(t)}=m^*;$$
(ii)
$$\zeta_4(t)A_*(\rho)\leq A_*(\rho t)\leq \zeta_5(t)A_*(\rho), \quad \mbox{ for all } \rho, t\geq 0;$$
(iii)
$$\zeta_4(\|u\|_{A_*})\leq\int_{\mathbb{R}^N}A_*(|u|)dx\leq\zeta_5(\|u\|_{A_*}), \quad \mbox{ for all }u \in L^{A_*}(\mathbb{R}^N),$$
where $A_*$ is the Sobolev conjugate function of $A$, which is defined by }
$$A_{*}^{-1}(t)=\int_{0}^{t}\frac{A^{-1}(s)}{s^{\frac{N+1}{N}}}ds,\quad \mbox{ \it for } t\geq 0.$$

\vskip2mm
\noindent
 {\bf Proposition 2.5.} (\cite{Adams2003})\ \  {\it Under the assumptions of Lemma 2.4, the embedding
 $$W^{1,A}(\mathbb{R}^N)\hookrightarrow L^B(\mathbb{R}^N)$$
 is continuous for any $N$-function $B$ satisfying
 $$\limsup_{r\rightarrow 0}\frac{B(r)}{A(r)}<\infty \quad \mbox{ and } \quad \limsup_{r\rightarrow \infty}\frac{B(r)}{A_*(r)}<\infty.$$
 Therefore, there exists a constant C such that}
 \begin{eqnarray*}\label{ 2.5+}
 \|u\|_B\leq C\|u\|_{1, A}, \quad \mbox{ \it for all } u \in W^{1, A}(\mathbb{R}^N).
 \end{eqnarray*}
 \vskip2mm
 \noindent
   {\bf Remark 2.1.}   By Lemma 2.1 and Lemma 2.3, $(\phi_1)$-$(\phi_2)$ imply that $\Phi_i$ and $\widetilde{\Phi}_i$, $i=1,2$, are $N$-functions that satisfy $\Delta_2$-condition globally. Thus, $L^{\Phi_i}(\mathbb{R}^N)$ and $W^{1, \Phi_i}(\mathbb{R}^N)$ are separable and reflexive Banach spaces (see \cite{Adams2003, M. M. Rao2002}). By (ii) in Lemma 2.1, (ii) in Lemma 2.4 and Proposition 2.5, it is easy to obtain the embedding
 \begin{eqnarray*}\label{ 2.4}
 W^{1, \Phi_i}(\mathbb{R}^N)\hookrightarrow L^{p_i}(\mathbb{R}^N)
 \end{eqnarray*}
  is continuous with $p_i\in [m_i,l_i^*]$, $i=1,2$, if we let $A(r): = \Phi_i(r)$ and $B(r):=r^{p_i}$. Hence, there exist positive constants  $C_{0,i}$,  $i=1,2$,   such that
  \begin{eqnarray}\label{2.5}
 \|u\|_{{p_i}}\le C_{0,i}\|u\|_{1,\Phi_i},
 \end{eqnarray}
whenever $p_i\in [m_i,l_i^*]$, $i=1,2$.

\par
Next, we recall a variant of mountain pass theorem. Let $X$ be a Banach space. $\varphi \in C^{1}(X,\R)$ and $c\in\R$. A sequence $\{u_{n}\}\subset X$ is called (PS)$_{c}$ sequence if the sequence $\{u_{n}\}$ satisfies
 \begin{equation*}
   \varphi(u_{n})\to c,\,\ \varphi'(u_{n})\to 0 .
   \end{equation*}
   \vskip2mm
 \noindent
 {\bf Lemma 2.5.} (Mountain Pass Theorem \cite{Brezis1991, Rabinowitz1986, Willem M1997}) Let X be a Banach space, $\varphi \in C^{1}(X,\R)$, $w\in X$ and $r>0$ be such that  $\|w\|>r$ and
 $$
 b:=\inf_{\|u\|=r} \varphi(u)>\varphi(0)\geq \varphi(w).
 $$
 Then there exists a (PS)$_{c}$ sequence with
  \begin{equation}\label{x1}
 c:=\inf_{\gamma\in\Gamma}\max_{t\in[0,1]}\varphi(\gamma(t)),
 \end{equation}
 and
 $$
 \Gamma:=\{\gamma\in([0,1],X):\gamma(0)=0,\gamma(1)=w\}.
 $$
 \vskip2mm

 \vskip2mm
 {\section{Proofs}}
  \setcounter{equation}{0}
  \par
  Define $W:=W^{1,\Phi_1}(\R^N)\times W^{1,\Phi_2}(\R^N)$ with the norm
  $$
  \|(u,v)\|= \|u\|_{1,\Phi_1}+ \|v\|_{1,\Phi_2}=\|\nabla u\|_{\Phi_1}+\|u\|_{\Phi_1}+\|\nabla v\|_{\Phi_2}+\|v\|_{\Phi_2}.
  $$
  Then $(W,\|\cdot\|)$ is a separable and reflexive Banach space.
  \par
  Next we use the idea in \cite{costa} to prove our theorem, on the whole,  which is the cut-off technique together with energy estimates. In order to adapt system (\ref{aa1}), we make an extension to $\R^2$ for the cut-off function in \cite{costa}.
  For some $\delta>0$, let $\rho_\delta\in C^1(\R\times \R,[0,1])$ be a cut-off function defined by
  \begin{eqnarray}\label{d2}
 \rho_\delta(t,s)= \begin{cases}
           1, \;\;\;  \text { if }\;\;|(t,s)| \le \delta/2, \\
           0, \;\;\;  \text { if }\;\;|(t,s)| \ge \delta
          \end{cases}
   \end{eqnarray}
  and $t\rho_t'(t,s)+s\rho_s'(t,s)\le 0 $ for all $(t,s)\in \R^2$. We  give some examples and their figures for the cut-off functions as follows
     \begin{eqnarray}
 \label{d3}(1)\ \ \rho_\delta(t,s)= \begin{cases}
           1, &  \text { if }\;\;|(t,s)|\le \frac{\delta}{2}, \\
          \sin\dfrac{8\pi(t^2+s^2-\delta^2)^2}{9\delta^4},    &  \text { if }\;\; \frac{\delta}{2}< |(t,s)|<\delta, \\
           0,&  \text { if }\;\;|(t,s)|\ge \delta;
          \end{cases}\\
 \label{d4}(2)\ \ \rho_\delta(t,s)= \begin{cases}
           1, &  \text { if }\;\;|(t,s)| \le \frac{\delta}{2}, \\
          \sin^2\dfrac{2\pi(t^2+s^2-\delta^2)}{3\delta^2},    & \text { if }\;\; \frac{\delta}{2}< |(t,s)|< \delta, \\
           0, & \text { if }\;\;|(t,s)|\ge \delta;
          \end{cases}\\
 \label{d5}(3)\ \ \rho_\delta(t,s)= \begin{cases}
           1, &  \text { if }\;\;|(t,s)| \le \frac{\delta}{2}, \\
          \cos\dfrac{8\pi(t^2+s^2-\frac{\delta^2}{4})^2}{9\delta^4},    &  \text { if }\;\; \frac{\delta}{2}< |(t,s)|< \delta, \\
           0,&  \text { if }\;\;|(t,s)|\ge\delta;
          \end{cases}\\
 \label{d6}(4)\ \ \rho_\delta(t,s)= \begin{cases}
           1, &  \text { if }\;\;|(t,s)| \le \frac{\delta}{2}, \\
          \cos^2\dfrac{2\pi(t^2+s^2-\frac{\delta^2}{4})}{3\delta^2},    & \text { if }\;\; \frac{\delta}{2}< |(t,s)|<\delta, \\
           0, & \text { if }\;\;|(t,s)|\ge \delta.
          \end{cases}
   \end{eqnarray}
 \begin{figure}[htbp]
\centering
\subfigure[Example (1) with $\delta$= 4.]{
\includegraphics[width=5.5cm]{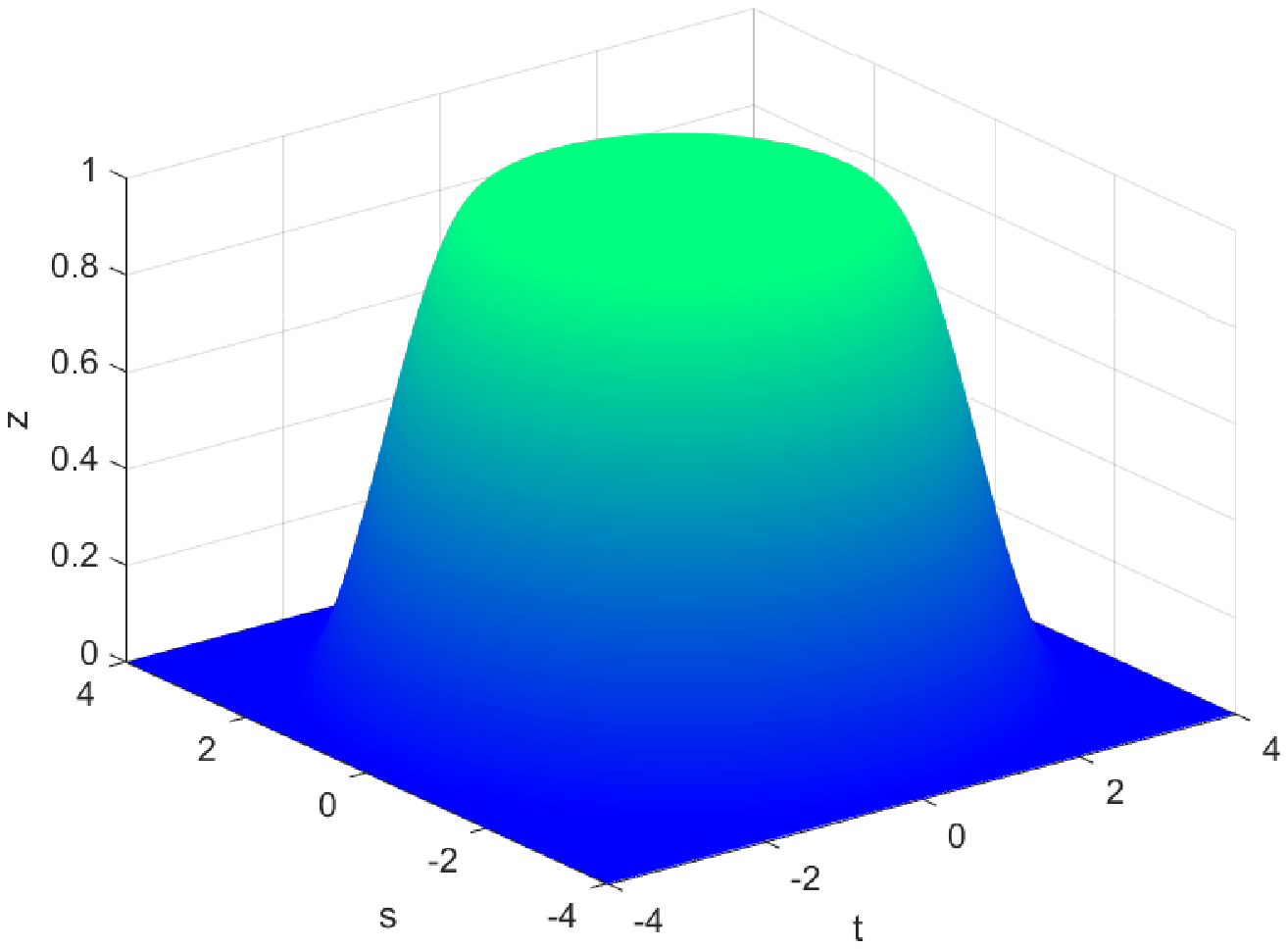}
}
\quad
\subfigure[Example (1) with $\delta$= 4.]{
\includegraphics[width=5.5cm]{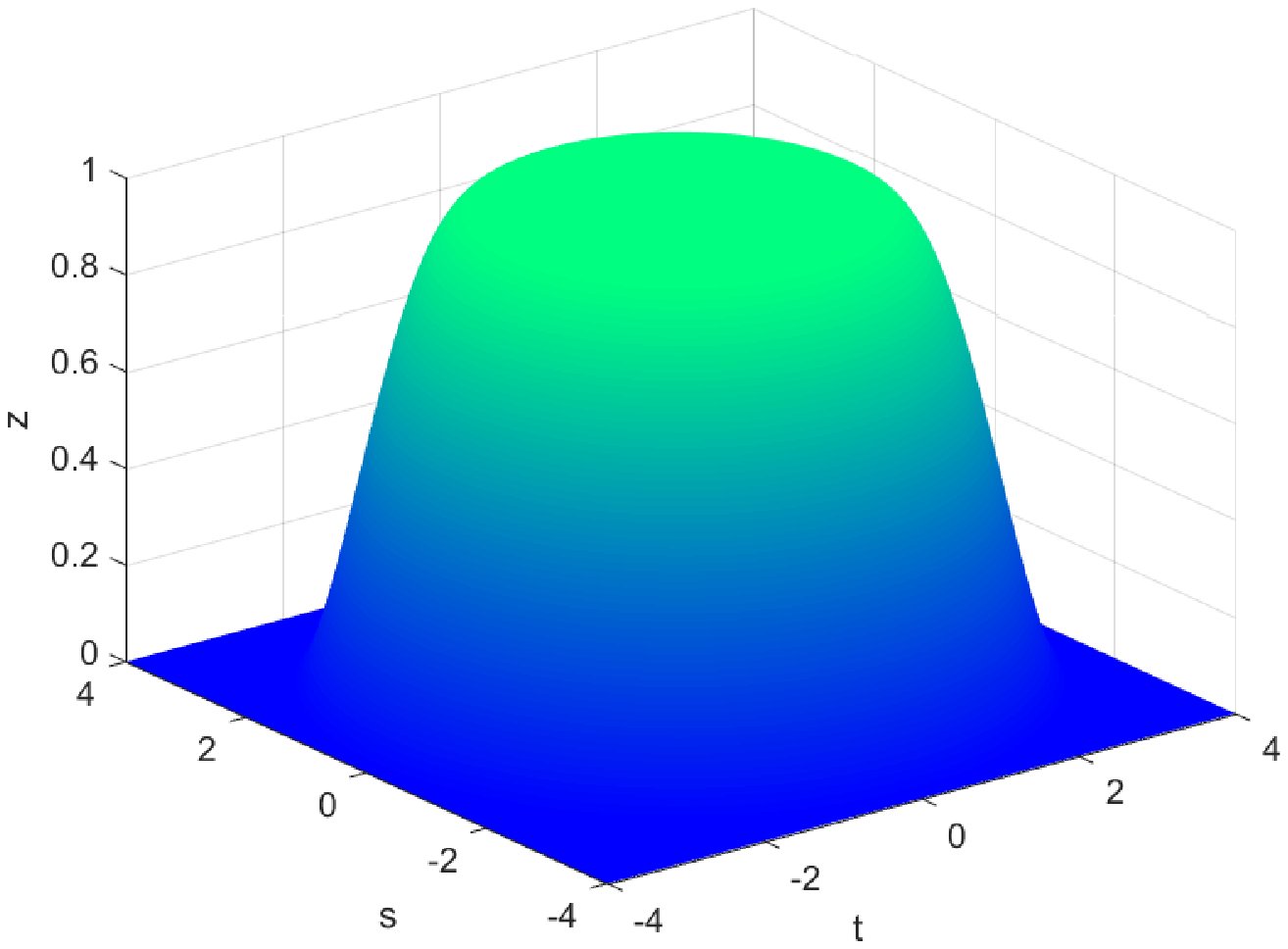}
}
\quad
\subfigure[Example (3) with $\delta$= 4.]{
\includegraphics[width=5.5cm]{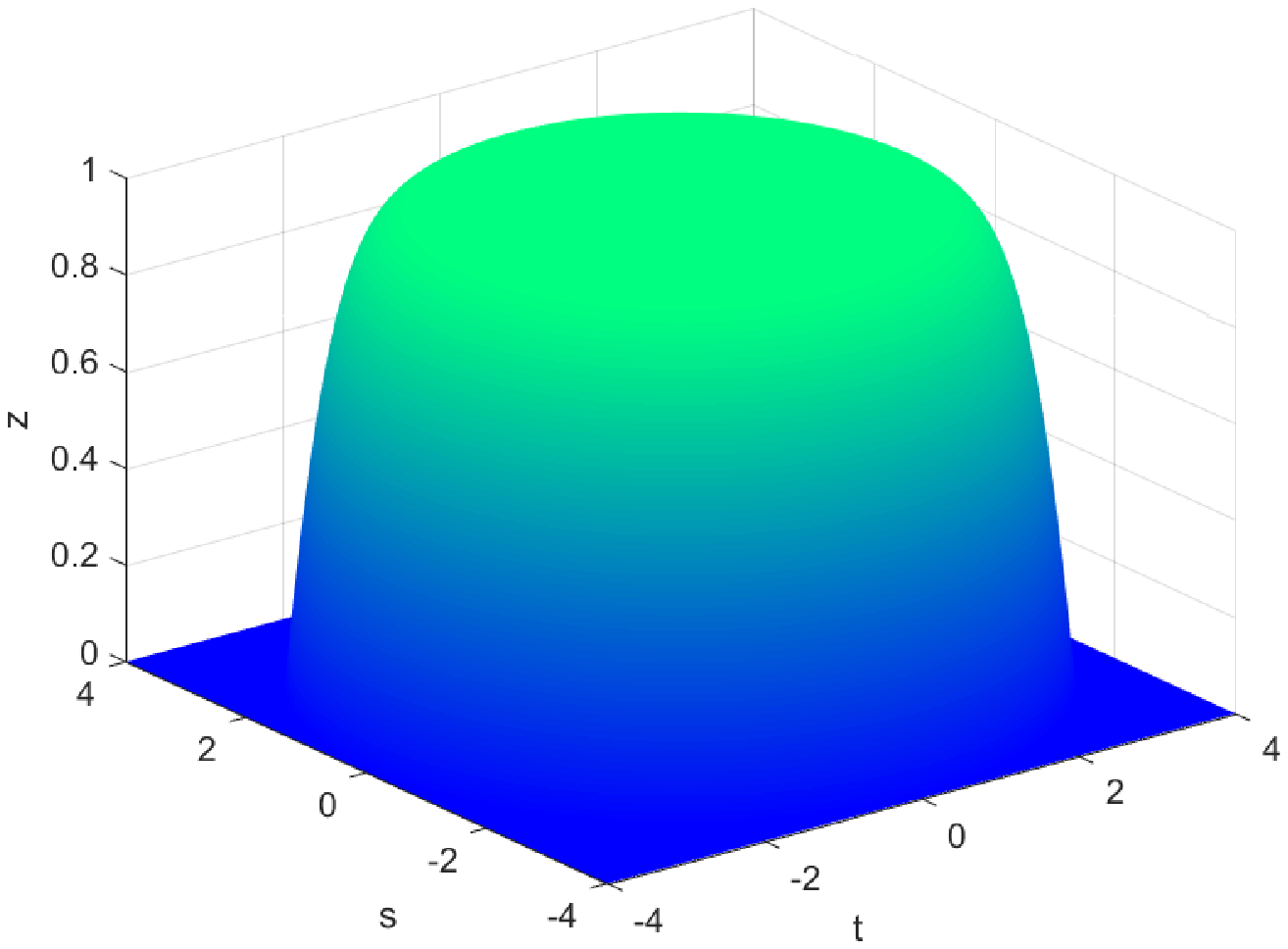}
}
\quad
\subfigure[Example (4) with $\delta$= 4.]{
\includegraphics[width=5.5cm]{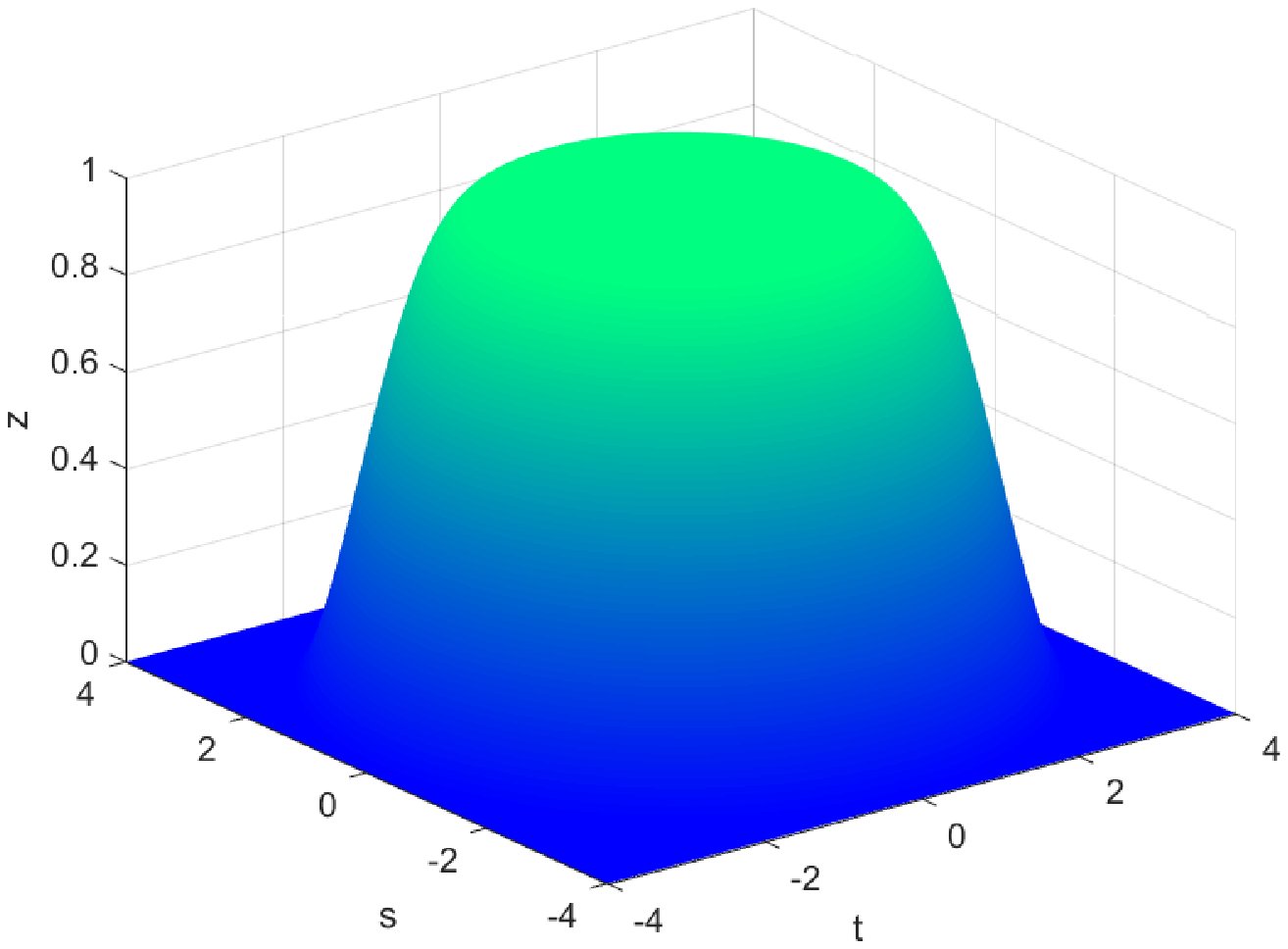}
}

\end{figure}
\vskip2mm
\noindent
 {\bf Remark 3.1.} For the examples of $\rho_\delta$, it seems to be natural to use  exponential functions as a link between $1$ and $0$ because of their infinite differentiability. However, we can not find such examples with exponential functions because it seems to be difficult to ensure that $\rho$ is differentiable at both $|(t,s)|=\frac{\delta}{2}$ and $|(t,s)|=\delta$. From the characteristic or shape of $\rho_\delta$, sine functions and cosine functions seem to be better choices.

 \vskip2mm
   By (F0)-(F3) and similar to the argument of Remark 2.8 in \cite{Wang2016}, there exist positive constants $M_5$ and $M_6$ such that
  \begin{eqnarray}\label{a8}
 |F(x,t,s)|\le M_5|t|^{r_1} +M_6|s|^{r_2}
   \end{eqnarray}
   for all $|(t,s)|< 2$ and $x\in  \R^N$.
   In fact, it follows from (F0) and (F3) that $F(x,0,0)=0$ for all $x\in \R^N$. Then if $r_2\ge r_1$, $t\in (-2,0)$ and $s\in (-\sqrt{4-t^2},0)$, by (F2), we have
    \begin{eqnarray}\label{cc1}
 |F(x,t,s)|
 &\leq &    \int_{t}^{0}|F_\tau(x,\tau,s)|d\tau+\int_{s}^{0}|F_\varsigma(x,0,\varsigma)|d\varsigma\nonumber\\
 &\leq &    \int_{t}^{0}\left(M_3 |\tau|^{r_1-1}+M_4 |s|^{r_2-1}\right)d\tau
           +\int_{s}^{0}M_4 |\varsigma|^{r_2-1}d\varsigma\nonumber\\
 &\leq  &   \frac{M_3}{r_1} |t|^{r_1}+M_4 |s|^{r_2-1}|t| +\frac{M_4}{r_2} |s|^{r_2}     \nonumber\\
 &\leq &   \frac{M_3}{r_1} |t|^{r_1}+\frac{M_4(r_2-1)}{r_2} |s|^{r_2}+\frac{M_4}{r_2}|t|^{r_2} +\frac{M_4}{r_2} |s|^{r_2} \nonumber\\
 &= &   \left(\frac{M_3}{r_1}+\frac{M_4}{r_2}|t|^{r_2-r_1}\right) |t|^{r_1} +M_4 |s|^{r_2}\nonumber\\
 &\leq &   \left(\frac{M_3}{r_1}+\frac{2^{r_2-r_1}M_4}{r_2}\right) |t|^{r_1} +M_4|s|^{r_2}
 \end{eqnarray}
 for all $x\in  \R^N$. If $r_1> r_2$, $t\in (-2,0)$ and $s\in (-\sqrt{4-t^2},0)$, we have
    \begin{eqnarray}\label{cc2}
 |F(x,t,s)|
 &\leq &    \int_{s}^{0}|F_\varsigma(x,t,\varsigma)|d\varsigma+\int_{t}^{0}|F_\tau(x,\tau,0)|d\tau\nonumber\\
 &\leq &    \int_{s}^{0}\left(M_3 |t|^{r_1-1}+M_4 |\varsigma|^{r_2-1}\right)d\varsigma
           +\int_{t}^{0}M_3 |\tau|^{r_1-1}d\tau\nonumber\\
 &\leq  &   M_3|t|^{r_1-1}|s|+\frac{M_4}{r_2} |s|^{r_2} +\frac{M_4}{r_1} |t|^{r_1}      \nonumber\\
 &\leq &   \frac{M_3(r_1-1)}{r_1} |t|^{r_1}+\frac{M_3}{r_1} |s|^{r_1}+\frac{M_4}{r_2}|s|^{r_2} +\frac{M_3}{r_1} |t|^{r_1}\nonumber\\
 &= &   \frac{M_3}{r_1}|t|^{r_1} + \left(\frac{M_4}{r_2}+\frac{M_3}{r_1}|s|^{r_1-r_2}\right)|s|^{r_2}\nonumber\\
 &\leq &  M_3|t|^{r_1} + \left(\frac{M_4}{r_2}+\frac{2^{r_1-r_2}M_3}{r_1}\right) |s|^{r_2}
 \end{eqnarray}
 for all $x\in  \R^N$. Combing (\ref{cc1}) and (\ref{cc2}), let $M_5=M_3+2^{|r_1-r_2|}M_4$ and  $M_6=M_4+2^{|r_1-r_2|}M_3$. Then (\ref{a8}) holds.
 Similar arguments can be done for other cases that $|(t,s)|<2$.
 \par
  Let $\delta=4$ in (\ref{d2}).
   Define $\widetilde{F}:\R^N\times\R\times\R\to \R$ by
$$
\widetilde{F}(x,t,s)=\rho_4(t,s)F(x,t,s)+(1-\rho_4(t,s)) M_5|t|^{r_1}+(1-\rho_4(t,s))M_6|s|^{r_2}.
$$
Then by  (F0)-(F3), the definition of $\rho_4$ and a direct computation, it is easy to obtain the following lemma:
 \vskip2mm
 \noindent
   {\bf Lemma 3.1.} {\it Assume that (F0)-(F3) hold. Then \\
    (F0)$'$ $\widetilde{F}\in C^1(\R^N\times \R\times \R,\R)$, $\widetilde{F}$ is $1-$periodic in $x\in \R^N$ and $\widetilde{F}(x,0,0)=0$ for all $x\in \R^N$;\\
  (F1)$'$
  \begin{eqnarray*}
 0\le\widetilde{F}(x,t,s)\leq M_5 |t|^{r_1}+M_6 |s|^{r_2},  \ \ \mbox{for all }(t,s)\in \mathbb {R}^2 \mbox{ and }x\in \R^N;
  \end{eqnarray*}
 (F2)$'$ there exist $M_7>0$ and $M_8>0$  such that
 \begin{eqnarray*}
 |\widetilde{F}_t(x,t,s)|\le M_7|t|^{r_1-1}+M_8|s|^{r_2-1},\ \  |\widetilde{F}_s(x,t,s)|\le  M_7|t|^{r_1-1}+M_8|s|^{r_2-1}
   \end{eqnarray*}
 for all $(t,s)\in \mathbb {R}^2$ and $x\in \R^N$;\\
 (F3)$'$
   \begin{eqnarray*}
 \theta \widetilde{F}(x,t,s)\leq \widetilde{F}_t(x,t,s)t+\widetilde{F}_s(x,t,s)s,  \ \ \mbox{for all }(t,s)\in \R^2/\{(0,0)\} \mbox{ and }x\in \R^N,
  \end{eqnarray*}
  where $\theta=\min\{r_1,r_2,\mu_1,\mu_2\}$.
 }

 \vskip2mm
 \par
   Consider the modified problem
   \begin{eqnarray}\label{aaaa1}
\begin{cases}
  -\mbox{div}(\phi_1(|\nabla u|)\nabla u)+V_1(x)\phi_1(|u|)u=\lambda  \widetilde{F}_u(x, u,v), \ \ x\in \R^N,\\
   -\mbox{div}(\phi_2(|\nabla v|)\nabla v)+V_2(x)\phi_2(|v|)v=\lambda  \widetilde{F}_v(x, u,v), \ \ x\in \R^N,\\
  u\in W^{1,\Phi_1}(\R^N), v\in W^{1,\Phi_2}(\R^N).
   \end{cases}
 \end{eqnarray}
Define the functional $\widetilde{J}_{\lambda}: W\to \R$ by
  $$
  \widetilde{J}_{\lambda}(u,v)=\int_{\R^N}\Phi_1(|\nabla u|)dx+\int_{\R^N}\Phi_2(|\nabla v|)dx+\int_{\R^N}V_1(x)\Phi_1(|u|)dx+\int_{\R^N}V_2(x)\Phi_2(|v|)dx-\lambda\int_{\R^N} \widetilde{F}(x,u,v)dx.
  $$
  By (F1)$'$ and a standard procedure, we can obtain that $ \widetilde{J}_{\lambda} $ is well defined and $\widetilde{J}_{\lambda}\in C^1(W,\R)$ and
\begin{eqnarray*}
          \langle \widetilde{J}'_{\lambda}(u,v),(\tilde{u},\tilde{v})\rangle
 &  = &  \int_{\R^N}(\phi_1(|\nabla u|)\nabla u,\nabla \tilde{u})dx+\int_{\R^N}V_1(x)\phi_1(|u|)u\tilde{u}dx\nonumber\\
  &  &   +\int_{\R^N}(\phi_2(|\nabla v|)\nabla v,\nabla \tilde{v})dx+\int_{\R^N}V_2(x)\phi_2(|v|)v\tilde{v}dx\nonumber\\
 &    &  -\lambda\int_{\R^N}\widetilde{F}_u(x,u,v)\tilde{u}dx -\lambda\int_{\R^N}\widetilde{F}_v(x,u,v)\tilde{v}dx
 \end{eqnarray*}
for all $(\tilde{u},\tilde{v})\in W$.
 \vskip2mm
 \noindent
   {\bf Lemma 3.2.} {\it $\widetilde{J}_{\lambda}$ satisfies the mountain pass geometry, that is,\\
   (i) there exist two positive constants $\gamma,\eta$ such that $\widetilde{J}_{\lambda}(u,v)\ge \eta$ for all $\|(u,v)\|=\gamma$; \\
   (ii) there exists $(u_0,v_0)\in C_0^\infty(\R^N)/\{0\}\times  C_0^\infty(\R^N)/\{0\}$ with $u_0>0$,  $v_0>0$ and $\|(u_0,v_0)\|_{\infty}:=\max_{x\in\R^N}(u_0^2+v_0^2)^{\frac{1}{2}}< 1$ such that $\widetilde{J}_{\lambda}(u_0,v_0)<0$.}
   \vskip2mm
 \noindent
 {\bf Proof.} If $|t|\le 1$, by (F2)$'$, $(\phi_2)$ and Lemma 2.1,  we have
   \begin{eqnarray}\label{b1}
         \frac{\widetilde{F}_t(x,t,s)}{\phi_1(|t|)|t|+\widetilde{\Phi_1}^{-1}(\Phi_2(|s|))}
  & \le & \left|\frac{\widetilde{F}_t(x,t,s)|t|}{\phi_1(|t|)t^2+|t|\widetilde{\Phi_1}^{-1}(\Phi_2(|s|))}\right|\nonumber\\
  &  \le & \frac{ M_7|t|^{r_1} +M_8|t||s|^{r_2-1}}{\phi_1(|t|)t^2}\nonumber\\
  &  \le & \frac{M_7|t|^{r_1} +M_8|t||s|^{r_2-1}}{l_1\Phi_1(|t|)} \nonumber\\
  &  \le & \frac{M_7|t|^{r_1} +M_8|t||s|^{r_2-1}}{l_1\min\{|t|^{l_1}, |t|^{m_1}\}\Phi(1)}\nonumber\\
  &   =  & \frac{M_7|t|^{r_1} +M_8|t||s|^{r_2-1}}{l_1\Phi(1)|t|^{m_1}}\nonumber\\
    &   =  & \frac{M_7|t|^{r_1-1} +M_8|s|^{r_2-1}}{l_1\Phi(1)|t|^{m_1-1}}.
   \end{eqnarray}
   Since $r_i>\max\{m_1,m_2\}, i=1,2$, (\ref{b1}) implies that
    \begin{eqnarray*}
       \lim_{|(t,s)|\to 0}  \frac{\widetilde{F}_t(x,t,s)}{\phi_1(|t|)|t|+\widetilde{\Phi}_1^{-1}(\Phi_2(|s|))}=0.
   \end{eqnarray*}
   Similarly, we also get that
    \begin{eqnarray*}
       \lim_{|(t,s)|\to 0}  \frac{\widetilde{F}_s(x,t,s)}{\widetilde{\Phi}_2^{-1}(\Phi_1(|t|))+\phi_2(|s|)|s|}=0.
   \end{eqnarray*}
   Next, we prove that
     \begin{eqnarray}\label{b2}
       \lim_{|(t,s)|\to \infty}  \frac{\widetilde{F}_t(x,t,s)}{\Phi_{1*}'(|t|)+\widetilde{\Phi}_{1*}^{-1}(\Phi_{2*}(|s|))}=0.
   \end{eqnarray}
   We divide three cases. For the case that $|t|\to \infty$ and $|s|$ is bounded, we assume that $|t|\ge 1$. Then by (F2)$'$, Lemma 2.4 and $l_1^*\le m_1^*$,  we have
   \begin{eqnarray}\label{c1}
           \frac{\widetilde{F}_t(x,t,s)}{\Phi_{1*}'(|t|)+\widetilde{\Phi}_{1*}^{-1}(\Phi_{2*}(|s|))}
   & \le &  \left|\frac{\widetilde{F}_{t}(x,t,s)|t|}{\Phi_{1*}'(|t|)|t|+|t|\widetilde{\Phi}_{1*}^{-1}(\Phi_{2*}(|s|))}\right|\nonumber\\
   &   \le & \frac{M_7|t|^{r_1} +M_8|t||s|^{r_2-1}}{\Phi_{1*}'(|t|)|t|}\nonumber\\
   &  \le & \frac{M_7|t|^{r_1} +M_8|t||s|^{r_2-1}}{l_1^*\min\{|t|^{l_1^*}, |t|^{m_1^*}\}\Phi_{1*}(1)}=\frac{M_7|t|^{r_1-1} +M_8|s|^{r_2-1}}{l_1^*\Phi_{1*}(1)|t|^{l_1^*-1}}.
   \end{eqnarray}
   Note that $l_1^*>r_1$. Then (\ref{b2}) holds. For the case that $|t|$ is bounded and $|s|\to \infty$, we assume that $|s|\ge 1$. Since $\frac{\widetilde{\Phi}_{1}^{-1}(s)}{s^{\frac{N+1}{N}}}$ is nondecreasing on $(0,\infty)$, we have $\frac{\widetilde{\Phi}_{1}^{-1}(s)}{s^{\frac{N+1}{N}}}\ge \widetilde{\Phi}_{1}^{-1}(1)$ for all $s\ge1$. Hence, by Lemma 2.4, we have
     \begin{eqnarray}\label{b3}
         \widetilde{\Phi}_{1*}^{-1}(\Phi_{2*}(|s|))
  & \ge & \widetilde{\Phi}_{1*}^{-1}(\Phi_{2*}(1)\min\{|s|^{l_2^*},|s|^{m_2^*}\})\nonumber\\
  & =   &  \widetilde{\Phi}_{1*}^{-1}(\Phi_{2*}(1)|s|^{l_2^*})\nonumber\\
  & =   & \int_0^{\Phi_{2*}(1)|s|^{l_2^*}}\frac{\widetilde{\Phi}_{1}^{-1}(s)}{s^{\frac{N+1}{N}}}ds\nonumber\\
   & \ge & \widetilde{\Phi}_{1}^{-1}(1)\Phi_{2*}(1)|s|^{l_2^*}.
   \end{eqnarray}
   Then  by (F2)$'$, Lemma 2.4, (\ref{b3}) and $l_1^*\le m_1^*$, similar to the argument of (\ref{c1}), we have
    \begin{eqnarray*}
           \frac{\widetilde{F}_t(x,t,s)}{\Phi_{1*}'(|t|)+\widetilde{\Phi}_{1*}^{-1}(\Phi_{2*}(|s|))}
   & \le & \frac{M_7|t|^{r_1-1} +M_8|s|^{r_2-1}}{ \Phi_{1}^{-1}(1)\Phi_{2*}(1)|s|^{l_2^*}}.
   \end{eqnarray*}
   Note that $l_2^*> r_2$. Then (\ref{b2}) holds. For the case that $|t|\to \infty$ and  $|s|\to \infty$,  we assume that  $|t|\ge 1$ and $|s|\ge 1$. Then
   \begin{eqnarray*}
           \frac{\widetilde{F}_t(x,t,s)}{\Phi_{1*}'(|t|)+\widetilde{\Phi}_{1*}^{-1}(\Phi_{2*}(|s|))}
   & \le &  \left|\frac{\widetilde{F}_{t}(x,t,s)|t|}{\Phi_{1*}'(|t|)|t|+|t|\widetilde{\Phi}_{1*}^{-1}(\Phi_{2*}(|s|))}\right|\nonumber\\
   &   \le & \frac{M_7|t|^{r_1} +M_8|t||s|^{r_2-1}}{\Phi_{1*}'(|t|)|t|+|t|\widetilde{\Phi}_{1*}^{-1}(\Phi_{2*}(|s|))}\nonumber\\
   &  \le & \frac{M_7|t|^{r_1} +M_8|t||s|^{r_2-1}}{l^*\min\{|t|^{l_1^*}, |t|^{m_1^*}\}\Phi_{1*}(1)+|t|\Phi_{1}^{-1}(1)\Phi_{2*}(1)|s|^{l_2^*}}\nonumber\\
   &   =  &\frac{M_7|t|^{r_1-1} +M_8|s|^{r_2-1}}{l_1^*\Phi_{1*}(1)|t|^{l_1^*-1}+\Phi_{1}^{-1}(1)\Phi_{2*}(1)|s|^{l_2^*}}.
   \end{eqnarray*}
   Note that  $l_1^*> r_1$ and $l_2^*> r_2$.  Then (\ref{b2}) holds. Hence,  (F2)$'$ implies (F2) of Lemma 3.14 in \cite{Wang2017}.  So the conclusion holds by  Lemma 3.14 and Lemma 3.15  in \cite{Wang2017}.\qed
 \vskip2mm
 \noindent
   {\bf Remark 3.1.}  There is a similar result of Lemma 3.2 in \cite{Wang2017} (see \cite{Wang2017}, Corollary 3.3) where the authors did not show the detail proof. Here, we present the proof for readers' convenience.
\vskip2mm
\par
 Since (F0)$'$-(F3)$'$ imply those conditions of Theorem 3.1 in  \cite{Wang2017},  it follows from Lemma 3.14-Lemma 3.16 in \cite{Wang2017} that system (\ref{aaaa1}) has a nontrivial solution $(u_{\lambda},v_{\lambda})$ such that $\widetilde{J}_{\lambda}(u_{\lambda},v_{\lambda})=c_\lambda$ with
 \begin{eqnarray*}\label{kkk1}
 c_{{\lambda}}:=\inf_{\gamma\in\Gamma}\max_{t\in[0,1]}\widetilde{J}_{\lambda}(\gamma(t)),
\end{eqnarray*}
and
 $$
 \Gamma:=\{\gamma\in([0,1],W):\gamma(0)=(0,0),\gamma(1)=(u_0,v_0)\}.
 $$

\vskip2mm
\par
Next, we  make an estimate for $\|u_{\lambda}\|_{1,\Phi_1}$ and $\|v_{\lambda}\|_{1,\Phi_2}$.
We introduce the functional $\bar{J}_{\lambda}: W\to \R$ as follows
\begin{eqnarray}\label{422}
         \bar{J}_{\lambda}(u)
 & = &  \int_{\R^N}\Phi_1(|\nabla u|)dx+\int_{\R^N}V_1(x)\Phi_1(|u|)dx+\int_{\R^N}\Phi_2(|\nabla v|)dx+\int_{\R^N}V_2(x)\Phi_2(|v|)dx \nonumber\\
 &   &  -\lambda M_1\int_{\R^N}|u|^{k_1} dx-\lambda M_2\int_{\R^N}|v|^{k_2} dx.
             \nonumber
 \end{eqnarray}
 \vskip2mm
 \noindent
{\bf Lemma 3.3.} There exist $\Lambda_0>0$ such that for each $\lambda>\Lambda_0$,  there exists $C_*>0$ such that
\begin{eqnarray*}
&  &  \|u_{\lambda}\|_{1,{\Phi_1}}\leq  C_*\max\left\{\lambda^{-\frac{1}{k_1-l_1}}+\lambda^{-\frac{l_2}{l_1(k_2-l_2)}},\lambda^{-\frac{l_1}{m_1(k_1-l_1)}}
+\lambda^{-\frac{l_2}{m_1(k_2-l_2)}}\right\},\\
&  &  \|v_{\lambda}\|_{2,{\Phi_2}}
 \le C_*\max\left\{\lambda^{-\frac{l_1}{l_2(k_1-l_1)}}+\lambda^{-\frac{1}{k_2-l_2}},\lambda^{-\frac{l_1}{m_2(k_1-l_1)}}
+\lambda^{-\frac{l_2}{m_2(k_2-l_2)}}\right\}.
 \end{eqnarray*}
\vskip2mm
\noindent
{\bf Proof.} By (V0) and (V1), we know that $V_i^\infty:=\max_{x\in \bar{\Omega}}V_i(x), i=1,2$ exist. Then by Lemma 2.2,  we have
\begin{eqnarray}\label{422}
         \bar{J}_{\lambda}(su_0,sv_0)
 & =&  \int_{\R^N}\Phi_1(s|\nabla u_0|)dx+\int_{\R^N}V_1(x)\Phi_1(s|u_0|)dx +\int_{\R^N}\Phi_2(s|\nabla v_0|)dx \nonumber\\
 &  &  +\int_{\R^N}V_2(x)\Phi_2(s|v_0|)dx -\lambda M_1s^{k_1}\int_{\R^N}|u_0|^{k_1} dx-\lambda M_2s^{k_2}\int_{\R^N}|v_0|^{k_2} dx
             \nonumber\\
 & \le & \max\{s^{l_1}\|\nabla u_0\|_{\Phi_1}^{l_1},s^{m_1}\|\nabla u_0\|_{\Phi_1}^{m_1}\}+V_1^\infty\max\{s^{l_1}\|u_0\|_{\Phi_1}^{l_1},s^{m_1}\|u_0\|_{\Phi_1}^{m_1}\}\nonumber\\
 &    &  + \max\{s^{l_2}\|\nabla v_0\|_{\Phi_2}^{l_2},s^{m_2}\|\nabla v_0\|_{\Phi_2}^{m_2}\}+V_2^\infty\max\{s^{l_2}\|v_0\|_{\Phi_2}^{l_2},s^{m_2}\|v_0\|_{\Phi_2}^{m_2}\}\nonumber\\
 &   &   -\lambda M_1s^{k_1}\int_{\R^N}|u_0|^{k_1} dx-\lambda M_2s^{k_2}\int_{\R^N}|v_0|^{k_2} dx  \nonumber\\
 & \le & s^{l_1}\left(\|\nabla u_0\|_{\Phi_1}^{l_1}+V_1^\infty \|u_0\|_{\Phi_1}^{l_1}+\|\nabla u_0\|_{\Phi_1}^{m_1}+V_1^\infty \|u_0\|_{\Phi_1}^{m_1}\right)\nonumber\\
 &   & +s^{l_2}\left(\|\nabla v_0\|_{\Phi_2}^{l_2}+V_2^\infty \|v_0\|_{\Phi_2}^{l_2}+\|\nabla v_0\|_{\Phi_2}^{m_2}+V_2^\infty \|v_0\|_{\Phi_2}^{m_2}\right)\nonumber\\
 &   & -\lambda M_1s^{k_1}\int_{\R^N}|u_0|^{k_1} dx-\lambda M_2s^{k_2}\int_{\R^N}|v_0|^{k_2} dx \nonumber
 \end{eqnarray}
for $s\in [0,1]$ and $(u_0,v_0)$ obtained in Lemma 3.2. Let
\begin{eqnarray*}
g_1(s)&= & s^{l_1}\left(\|\nabla u_0\|_{\Phi_1}^{l_1}+V_1^\infty \|u_0\|_{\Phi_1}^{l_1}+\|\nabla u_0\|_{\Phi_1}^{m_1}+V_1^\infty \|u_0\|_{\Phi_1}^{m_1}\right)-\lambda M_1s^{k_1}\int_{\R^N}|u_0|^{k_1} dx,\nonumber\\
g_2(s)&= & s^{l_2}\left(\|\nabla v_0\|_{\Phi_2}^{l_2}+V_2^\infty \|v_0\|_{\Phi_2}^{l_2}+\|\nabla v_0\|_{\Phi_2}^{m_2}+V_2^\infty \|v_0\|_{\Phi_2}^{m_2}\right)-\lambda M_2s^{k_2}\int_{\R^N}|v_0|^{k_2} dx.
 \end{eqnarray*}
Obviously, there exist $\Lambda_1>0$ and $\Lambda_2>0$ such that $s_{\Lambda_1}:=\left(\frac{l_1\left(\|\nabla u_0\|_{\Phi_1}^{l_1}+V_1^\infty \|u_0\|_{\Phi_1}^{l_1}+\|\nabla u_0\|_{\Phi_1}^{m_1}+V_2^\infty \|u_0\|_{\Phi_1}^{m_1}\right)}{\lambda M_1k_1\int_{\R^N}|u_0|^{k_1}  dx}\right)^{\frac{1}{k_1-l_1}}\in (0,1)$ and  $s_{\Lambda_2}:=\left(\frac{l_2\left(\|\nabla v_0\|_{\Phi_2}^{l_2}+V_2^\infty \|v_0\|_{\Phi_2}^{l_2}+\|\nabla v_0\|_{\Phi_2}^{m_2}+V_2^\infty \|v_0\|_{\Phi_2}^{m_2}\right)}{\lambda M_2k_2\int_{\R^N}|v_0|^{k_2}  dx}\right)^{\frac{1}{k_2-l_2}}\in (0,1)$ if $\lambda>\max\{\Lambda_1,\Lambda_2\}$. Then when $s=s_{\Lambda_1}$,  $g_1(s)$ attains the maximum on $[0,1]$. Then there exists $C_{1,*}>0$ such that
$$
     \max_{s\in[0,1]} g_1(s)
=  \lambda^{-\frac{l_1}{k_1-l_1}}k_1^{-\frac{k_1}{k_1-l_1}}l_1^{\frac{l_1}{k_1-l_1}}\left(k_1-l_1\right)A^{\frac{k_1}{k_1-l_1}}
\left(M_1\int_{\R^N}|u_0|^{k_1}  dx\right)^{-\frac{l_1}{k_1-l_1}}
\le
C_{1,*}\lambda^{-\frac{l_1}{k_1-l_1}},
$$
where $A=\|\nabla u_0\|_{\Phi_1}^{l_1}+V_1^\infty \|u_0\|_{\Phi_1}^{l_1}+\|\nabla u_0\|_{\Phi_1}^{m_1}+V_1^\infty \|u_0\|_{\Phi_1}^{m_1}$.
When $s=s_{\Lambda_2}$,  $g_2(s)$ attains the maximum on $[0,1]$ and there exists $C_{2,*}>0$ such that
$$
\max_{s\in[0,1]} g_2(s)=\lambda^{-\frac{l_2}{k_2-l_2}}k_2^{-\frac{k_2}{k_2-l_2}}l_2^{\frac{l_2}{k_2-l_2}}\left(k_2-l_2\right)B^{\frac{k_2}{k_2-l_2}} \left(M_2\int_{\R^N}|v_0|^{k_2}  dx\right)^{-\frac{l_2}{k_2-l_2}}\le C_{2,*}\lambda^{-\frac{l_2}{k_2-l_2}},
$$
where $B=\|\nabla v_0\|_{\Phi_2}^{l_2}+V_2^\infty \|v_0\|_{\Phi_2}^{l_2}+\|\nabla v_0\|_{\Phi_2}^{m_2}+V_2^\infty \|v_0\|_{\Phi_2}^{m_2}.$
Notice that $|s(u_0(x),v_0(x))|\le \|(u_0,v_0)\|_{\infty}< 1$ for all $x\in \R^N$ and $s\in [0,1]$. Then by the definition of $\widetilde{F}$ and (F1), we have $\widetilde{F}(su_0,sv_0)={F}(su_0,sv_0)\ge M_1|su_0(x)|^{k_1}+M_2|sv_0(x)|^{k_2}$. Thus,
  \begin{eqnarray}\label{a1}
 c_{{\lambda}} :=  \inf_{\gamma\in\Gamma}\max_{s\in[0,1]}\widetilde{J}_{\lambda}(\gamma(s))\le \max_{s\in[0,1]}\widetilde{J}_{\lambda}(su_0,sv_0)\le \max_{s\in[0,1]}  \bar{J}_{\lambda}(su_0,sv_0)
\le C_{1,*}\lambda^{-\frac{l_1}{k_1-l_1}}+ C_{2,*}\lambda^{-\frac{l_2}{k_2-l_2}}.
\end{eqnarray}
Note that $(u_{\lambda},v_{\lambda})$ is  a critical point of $ \widetilde{J}_{\lambda}$ with critical value
$c_{\lambda}$.
 Since $\langle \widetilde{J}'_{\lambda}(u_{\lambda},v_{\lambda}),(u_{\lambda},v_{\lambda})\rangle =0$ and $\theta\ge \max\{m_1,m_2\}$, by  $(\phi_2)$, (F3)$'$ and Lemma 2.2, we have
 \begin{eqnarray}\label{k6}
 \theta c_{\lambda} &= &     \theta  \widetilde{J}_{\lambda}(u_{\lambda},v_{\lambda})\nonumber \\
 &= & \theta   \langle \widetilde{J}_{\lambda}(u_{\lambda},v_{\lambda}),(u_{\lambda},v_{\lambda})\rangle
         -\langle \widetilde{J}_{\lambda}^\prime(u_{\lambda},v_{\lambda}),(u_{\lambda},v_{\lambda})\rangle
          \nonumber \\
 & = &  \theta  \int_{\R^N}\Phi_1(|\nabla u_{\lambda}|)dx+\theta  \int_{\R^N}V_1(x)\Phi_1(|u_{\lambda}|)dx+\theta  \int_{\R^N}\Phi_2(|\nabla v_{\lambda}|)dx+\theta  \int_{\R^N}V_2(x)\Phi_2(|v_{\lambda}|)dx
          \nonumber \\
 &  &
          -\theta  \lambda\int_{\R^N} \widetilde{F}(x,u_{\lambda},v_{\lambda})dx -\int_{\R^N}\phi_1(|\nabla u_{\lambda}|)|\nabla u_{\lambda}|^2dx-\int_{\R^N}V_1(x)\phi_1(|u_{\lambda}|)u_{\lambda}^2dx\nonumber \\
 &   &  -\int_{\R^N}\phi_2(|\nabla v_{\lambda}|)|\nabla v_{\lambda}|^2dx-\int_{\R^N}V_2(x)\phi_2(|v_{\lambda}|)v_{\lambda}^2dx+\lambda\int_{\R^N} \widetilde{F}_u(x,u_{\lambda},v_{\lambda})u_{\lambda}dx\nonumber\\
 &   & +\lambda\int_{\R^N} \widetilde{F}_v(x,u_{\lambda},v_{\lambda})v_{\lambda}dx
              \nonumber \\
 &\geq & (\theta-m_1)  \int_{\R^N}\Phi_1(|\nabla u_{\lambda}|)dx+(\theta-m_1)   \int_{\R^N}V_1(x)\Phi_1(|u_{\lambda}|)dx+ (\theta-m_2)  \int_{\R^N}\Phi_2(|\nabla v_{\lambda}|)dx \nonumber \\
 &     &+(\theta-m_2)   \int_{\R^N}V_2(x)\Phi_2(|v_{\lambda}|)dx
             \nonumber \\
 &\geq & (\theta-m_1)\min\{\|\nabla u_{\lambda}\|_{\Phi_{1}}^{l_1},\|\nabla u_{\lambda}\|_{\Phi_1}^{m_1}\}+(\theta-m_1)V_{1,\infty}\min\{\|u_{\lambda}\|_{\Phi_1}^{l_1},\|u_{\lambda}\|_{\Phi_1}^{m_1}\}\nonumber\\
 &     & +(\theta-m_2)\min\{\|\nabla v_{\lambda}\|_{\Phi_2}^{l_2},\|\nabla v_{\lambda}\|_{\Phi_1}^{m_2}\}+(\theta-m_2)V_{2,\infty}\min\{\|v_{\lambda}\|_{\Phi_2}^{l_2},\|v_{\lambda}\|_{\Phi_2}^{m_2}\},
             \nonumber
 \end{eqnarray}
 which, together with (\ref{a1}), implies that
 \begin{eqnarray*}\label{k6}
 &    &  \min\{\|\nabla u_{\lambda}\|_{\Phi_1}^{l_1},\|\nabla u_{\lambda}\|_{\Phi_1}^{m_1}\}+V_{1,\infty}\min\{\|u_{\lambda}\|_{\Phi_1}^{l_1},\|u_{\lambda}\|_{\Phi_1}^{m_1}\} \nonumber\\
 & \le &  \frac{\theta}{\theta-m_1} c_{\lambda} \le  \frac{\theta}{\theta-m_1}C_{1,*}\lambda^{-\frac{l_1}{k_1-l_1}}+ \frac{\theta}{\theta-m_1}C_{2,*}\lambda^{-\frac{l_2}{k_2-l_2}}         \nonumber
  \end{eqnarray*}
 and
   \begin{eqnarray*}\label{k6}
 &    &  \min\{\|\nabla v_{\lambda}\|_{\Phi_2}^{l_2},\|\nabla v_{\lambda}\|_{\Phi_2}^{m_2}\}+V_{2,\infty}\min\{\|v_{\lambda}\|_{\Phi_2}^{l_2},\|v_{\lambda}\|_{\Phi_2}^{m_2}\} \nonumber\\
 & \le &  \frac{\theta}{\theta-m_2} c_{\lambda} \le   \frac{\theta}{\theta-m_2}C_{1,*}\lambda^{-\frac{l_1}{k_1-l_1}}+ \frac{\theta}{\theta-m_2}C_{2,*}\lambda^{-\frac{l_2}{k_2-l_2}}.
 \end{eqnarray*}
 Hence, there exists $C_*>0$ such that
\begin{eqnarray*}
          \|u_{\lambda}\|_{1,{\Phi_1}}
 &  =  &  \|\nabla u_{\lambda}\|_{\Phi_1}+\|u_{\lambda}\|_{\Phi_1}\\
 & \le & \left(\frac{1}{V_{1,\infty}}+1\right)\max\left\{\left(\frac{\theta}{\theta-m_1}C_{1,*}\right)^{\frac{1}{l_1}}\lambda^{-\frac{1}{k_1-l_1}}+\left(\frac{\theta}{\theta-m_1}C_{2,*}\right)^{\frac{1}{l_1}}\lambda^{-\frac{l_2}{l_1(k_2-l_2)}},\right.\\
&  &   \left.\left(\frac{\theta}{\theta-m_1}C_{1,*}\right)^{\frac{1}{m_1}}\lambda^{-\frac{l_1}{m_1(k_1-l_1)}}
+\left(\frac{\theta}{\theta-m_1}C_{2,*}\right)^{\frac{1}{m_1}}\lambda^{-\frac{l_2}{m_1(k_2-l_2)}}\right\}\\
&\le & C_*\max\left\{\lambda^{-\frac{1}{k_1-l_1}}+\lambda^{-\frac{l_2}{l_1(k_2-l_2)}},\lambda^{-\frac{l_1}{m_1(k_1-l_1)}}
+\lambda^{-\frac{l_2}{m_1(k_2-l_2)}}\right\}
 \end{eqnarray*}
 and
 \begin{eqnarray*}
          \|v_{\lambda}\|_{2,{\Phi_2}}
 &  =  &  \|\nabla v_{\lambda}\|_{\Phi_2}+\|v_{\lambda}\|_{\Phi_2}\\
 & \le & \left(\frac{1}{V_{2,\infty}}+1\right)\max\left\{\left(\frac{\theta}{\theta-m_2}C_{1,*}\right)^{\frac{1}{l_2}}\lambda^{-\frac{l_1}{l_2(k_1-l_1)}}
 +\left(\frac{\theta}{\theta-m_2}C_{2,*}\right)^{\frac{1}{l_2}}\lambda^{-\frac{1}{k_2-l_2}}\right.\\
&  &   \left.\left(\frac{\theta}{\theta-m_2}C_{1,*}\right)^{\frac{1}{m_2}}\lambda^{-\frac{l_1}{m_2(k_1-l_1)}}
+\left(\frac{\theta}{\theta-m_2}C_{2,*}\right)^{\frac{1}{m_2}}\lambda^{-\frac{l_2}{m_2(k_2-l_2)}}\right\}\\
& \le & C_*\max\left\{\lambda^{-\frac{l_1}{l_2(k_1-l_1)}}+\lambda^{-\frac{1}{k_2-l_2}},\lambda^{-\frac{l_1}{m_2(k_1-l_1)}}
+\lambda^{-\frac{l_2}{m_2(k_2-l_2)}}\right\}.
 \end{eqnarray*}
Let $\Lambda_0=\max\{\Lambda_1,\Lambda_2\}$. Then we complete the proof. \qed
 \vskip2mm
 \noindent
{\bf Lemma 3.4.} {\it There exists $\Lambda_*>\Lambda_0$ such that for all $\lambda>\Lambda_*$, }
\begin{eqnarray*}
\|u_{\lambda}\|_\infty\le 1, \|v_{\lambda}\|_\infty\le  1.
\end{eqnarray*}
{\bf Proof.} The proof is motivated by Lemma 2.6 in \cite{Medeiros}, which originates from \cite{Moser} and is often called Moser iteration technique (see also \cite{Chabrowski} and \cite{Figueiredo} for some related arguments). However, because of the interaction of $u$ and $v$, our proof is more difficult than those in \cite{Medeiros}, \cite{Chabrowski} and \cite{Figueiredo}, where some scalar equations were investigated.
 \par
 Next, we start the proof. Since $(u_{\lambda},v_{\lambda})$ is a critical point of $\widetilde{J}_{\lambda}$, we have
\begin{eqnarray}\label{kk1}
&  &  \int_{{\R^N}}(\phi_1(|\nabla u_{\lambda}|)\nabla u_{\lambda},\nabla \tilde{u})dx+ \int_{{\R^N}}V_1(x)\phi_1(|u_{\lambda}|)u_{\lambda}\tilde{u}dx=\lambda \int_{{\R^N}} \widetilde{F}_{u_{\lambda}}(x,u_{\lambda},v_{\lambda})\tilde{u}dx,\\
&  &   \int_{{\R^N}}(\phi_2(|\nabla v_{\lambda}|)\nabla v_{\lambda},\nabla \tilde{v})dx+ \int_{{\R^N}}V_2(x)\phi_2(|v_{\lambda}|)v_{\lambda}\tilde{v}dx=\lambda \int_{{\R^N}} \widetilde{F}_{v_{\lambda}}(x,u_{\lambda},v_{\lambda})\tilde{v}dx.\nonumber
\end{eqnarray}
for all $(\tilde{u},\tilde{v})\in W$. Next, we prove $\|u_{\lambda}\|_\infty\le 1$.    Without loss of generality, for each $k>0$, we define
$$
u_k=\begin{cases}
  u_{\lambda},  & \mbox{if } u_{\lambda} \le k,\\
  k,  & \mbox{if } u_{\lambda}  > k,
\end{cases}
$$
$\varphi_k=|u_k|^{l_1(\beta_1-1)}u_\lambda$ and $w_k=u_{\lambda}|u_k|^{\beta_1-1}$ with $\beta_1>1$.
By $(\phi_3)$, we have $\phi_1(|\nabla u_{\lambda}|)|\nabla u_{\lambda}|^2\ge q_{1} |\nabla u_{\lambda}|^{l_1}$. Then it follows from (\ref{u1}), (\ref{2.2tj}), H$\ddot{\mbox{o}}$lder inequality, Lemma 2.2 and (\ref{2.5}) that
\begin{eqnarray}\label{kk2}
&     &      q_{1} \int_{{\R^N}}|\nabla u_{\lambda}|^{l_1} |u_k|^{l_1(\beta_1-1)}dx\nonumber\\
& \le &  \int_{{\R^N}}\phi_1(|\nabla u_{\lambda}|)(\nabla u_{\lambda},\nabla \varphi_k)dx -l_1(\beta_1-1)
        \int_{{\R^N}}|u_k|^{l_1(\beta_1-1)-2}u_k u_\lambda\phi_1(|\nabla u_{\lambda}|)(\nabla u_{\lambda},\nabla u_k)dx\nonumber\\
& \le & - \int_{{\R^N}}V_1(x)\phi(|u_{\lambda}|)u_{\lambda}\varphi_kdx
       +\lambda \int_{{\R^N}} \widetilde{F}_{u_{\lambda}}(x,u_{\lambda},v_{\lambda})\varphi_kdx\nonumber\\
& \le &   \lambda M_7 \int_{{\R^N}} |u_{\lambda}|^{r_1}|u_k|^{l_1(\beta_1-1)}dx
        + \lambda M_8 \int_{{\R^N}} |v_{\lambda}|^{r_2-1}|u_\lambda||u_k|^{l_1(\beta_1-1)}dx\nonumber\\
& =  &   \lambda M_7 \int_{{\R^N}} |u_{\lambda}|^{r_1-l_1}|w_k|^{l_1}dx
      +  \lambda M_8 \int_{{\R^N}} |v_{\lambda}|^{r_2-1}|w_k||u_k|^{(l_1-1)(\beta_1-1)}dx\nonumber\\
& \le &  \lambda M_7 \int_{{\R^N}} |u_{\lambda}|^{r_1-l_1}|w_k|^{l_1}dx
      +  \lambda M_8 \int_{{\R^N}} |v_{\lambda}|^{r_2-1}|w_k||u_\lambda|^{(l_1-1)(\beta_1-1)}dx\nonumber\\
& =  &   \lambda M_7 \int_{{\R^N}} |u_{\lambda}|^{r_1-l_1}|w_k|^{l_1}dx
      +  \lambda M_8 \int_{{\R^N}} |v_{\lambda}|^{r_2-1}|u_\lambda|^{(l_1-1)(\beta_1-1)-\frac{r_1-l_1}{l_1}}|u_\lambda|^{\frac{r_1-l_1}{l_1}}|w_k|dx\nonumber\\
&\le  &\lambda M_7 \int_{{\R^N}} |u_{\lambda}|^{r_1-l_1}|w_k|^{l_1}dx
+\lambda M_8 \left(\int_{{\R^N}} |v_{\lambda}|^{\frac{l_1(r_2-1)}{l_1-1}}|u_{\lambda}|^{l_1(\beta_1-1)-\frac{r_1-l_1}{l_1-1}}dx\right)^{1-\frac{1}{l_1}}  \left(\int_{{\R^N}} |w_k|^{l_1}|u_\lambda|^{r_1-l_1}dx\right)^{\frac{1}{l_1}}\nonumber\\
&\le  &
\lambda M_8 \left(\int_{{\R^N}} |v_{\lambda}|^{\frac{\Theta_1 l_1(r_2-1)}{(\Theta_1-1)(l_1-1)}}dx\right)^{\frac{(l_1-1)(\Theta_1-1)}{\Theta_1 l_1}}
\left(\int_{{\R^N}} |u_{\lambda}|^{\Theta_1 l_1(\beta_1-1)-\frac{\Theta_1 (r_1-l_1)}{l_1-1}}dx\right)^{\frac{l_1-1}{\Theta_1 l_1}}
\left(\int_{{\R^N}} |w_k|^{l_1}|u_\lambda|^{r_1-l_1}dx\right)^{\frac{1}{l_1}}\nonumber\\
&& +\lambda M_7 \int_{{\R^N}} |u_{\lambda}|^{r_1-l_1}|w_k|^{l_1}dx\nonumber\\
&=  &
\lambda M_8 \left(\int_{{\R^N}} |v_{\lambda}|^{\frac{\Theta_1 l_1(r_2-1)}{(\Theta_1-1)(l_1-1)}-m_2}|v_{\lambda}|^{m_2}dx
            \right)^{\frac{(l_1-1)(\Theta_1-1)}{\Theta_1 l_1}}
           \left(\int_{{\R^N}} |u_{\lambda}|^{\Theta_1 l_1(\beta_1-1)-\frac{\Theta_1 (r_1-l_1)}{l_1-1}-m_1}|u_{\lambda}|^{m_1}dx
          \right)^{\frac{l_1-1}{\Theta_1 l_1}} \nonumber\\
&&   \times
       \left(\int_{{\Omega}} |w_k|^{l_1}|u_\lambda|^{r_1-l_1}dx\right)^{\frac{1}{l_1}}
     + \lambda M_7 \int_{{\R^N}} |u_{\lambda}|^{r_1-l_1}|w_k|^{l_1}dx\nonumber\\
&\le  &
   \lambda M_8
   \left(\int_{{\R^N}} |w_k|^{l_1}|u_\lambda|^{r_1-l_1}dx\right)^{\frac{1}{l_1}}
   \left(\int_{{\R^N}}
         \left(\Phi_{2}(|v_{\lambda}|)+\Phi_{2}^{-1}(1)\right)^{\frac{\Theta_1 l_1(r_2-1)}{(\Theta_1-1)(l_1-1)}-m_2}|v_{\lambda}|^{m_2}
     dx\right)^{\frac{(l_1-1)(\Theta_1-1)}{\Theta_1 l_1}}
\nonumber\\
&&  \times
     \left(\int_{{\R^N}}
          \left(\Phi_{1}(|u_{\lambda}|)+\Phi_{1}^{-1}(1)\right)^{\Theta_1 l_1(\beta_1-1)-\frac{\Theta_1 (r_1-l_1)}{l_1-1}-m_1}|u_{\lambda}|^{m_1}
       dx\right)^{\frac{l_1-1}{\Theta_1 l_1}}
   +\lambda M_7 \int_{{\R^N}} |u_{\lambda}|^{r_1-l_1}|w_k|^{l_1}dx\nonumber\\
&\le  &
   \lambda M_8 2^{r_2-1-\frac{(m_2+1)(l_1-1)(\Theta_1-1)}{\Theta_1 l_1}} 2^{(l_1-1)(\beta_1-1)-\frac{r_1-l_1}{l_1}-\frac{(m_1+1)(l_1-1)}{\Theta_1 l_1}}
   \left(\int_{{\R^N}} |w_k|^{l_1}|u_\lambda|^{r_1-l_1}dx\right)^{\frac{1}{l_1}}
\nonumber\\
&&    \times
   \left(\int_{{\R^N}}
          \left( \Phi_{2}^{\frac{\Theta_1 l_1(r_2-1)}{(\Theta_1-1)(l_1-1)}-m_2}(|v_{\lambda}|)
                +(\Phi_{2}^{-1}(1))^{\frac{\Theta_1 l_1(r_2-1)}{(\Theta_1-1)(l_1-1)}-m_2}
          \right)
       |v_{\lambda}|^{m_2}
     dx\right)^{\frac{(l_1-1)(\Theta_1-1)}{\Theta_1 l_1}}
\nonumber\\
&&
   \times
   \left(\int_{{\R^N}}
         \left( \Phi_{1}^{\Theta_1 l_1(\beta_1-1)-\frac{\Theta_1 (r_1-l_1)}{l_1-1}-m_1}(|u_{\lambda}|)
              +(\Phi_{1}^{-1}(1))^{\Theta_1 l_1(\beta_1-1)-\frac{\Theta_1 (r_1-l_1)}{l_1-1}-m_1}
         \right)
         |u_{\lambda}|^{m_1}
     dx\right)^{\frac{l_1-1}{\Theta_1 l_1}}
\nonumber\\
&&
   +\lambda M_7 \int_{{\Omega}} |u_{\lambda}|^{r_1-l_1}|w_k|^{l_1}dx\nonumber\\
& <  &
  \lambda M_8 2^{r_2-1-\frac{(m_2+1)(l_1-1)(\Theta_1-1)}{\Theta_1 l_1}} 2^{(l_1-1)(\beta_1-1)-\frac{r_1-l_1}{l_1}-\frac{(m_1+1)(l_1-1)}{\Theta_1 l_1}}
   \left(\int_{{\R^N}} |w_k|^{l_1}|u_\lambda|^{r_1-l_1}dx\right)^{\frac{1}{l_1}}
\nonumber\\
&&     \times
   \left(
      \int_{{\R^N}}
         \left(\Phi_{2}(|v_{\lambda}|) |v_{\lambda}|^{m_2}
         +(\Phi_{2}^{-1}(1))^{\frac{\Theta_1 l_1(r_2-1)}{(\Theta_1-1)(l_1-1)}-m_2} |v_{\lambda}|^{m_2}\right)
     dx\right)^{\frac{(l_1-1)(\Theta_1-1)}{\Theta_1 l_1}}
\nonumber\\
&&   \times
   \left(
      \int_{{\R^N}}
         \left(\Phi_{1}(|u_{\lambda}|)|u_{\lambda}|^{m_1}
         +(\Phi_{1}^{-1}(1))^{\Theta_1 l_1(\beta_1-1)-\frac{\Theta_1 (r_1-l_1)}{l_1-1}-m_1}|u_{\lambda}|^{m_1}\right)
     dx\right)^{\frac{l_1-1}{\Theta_1 l_1}}
\nonumber\\
&&
   +\lambda M_7 \int_{{\Omega}} |u_{\lambda}|^{r_1-l_1}|w_k|^{l_1}dx\nonumber\\
& =  &
   \lambda M_8 2^{r_2-1-\frac{(m_2+1)(l_1-1)(\Theta_1-1)}{\Theta_1 l_1}} 2^{(l_1-1)(\beta_1-1)-\frac{r_1-l_1}{l_1}-\frac{(m_1+1)(l_1-1)}{\Theta_1 l_1}}
   \left(\int_{{\R^N}} |w_k|^{l_1}|u_\lambda|^{r_1-l_1}dx\right)^{\frac{1}{l_1}}
\nonumber\\
&&  \times
    \left(
      \int_{{\R^N}}
           \Phi_{2}(|v_{\lambda}|) |v_{\lambda}|^{m_2}
      dx
   +
      \int_{{\R^N}}
(\Phi_{2}^{-1}(1))^{\frac{\Theta_1 l_1(r_2-1)}{(\Theta_1-1)(l_1-1)}-m_2} |v_{\lambda}|^{m_2}
      dx
     \right)^{\frac{(l_1-1)(\Theta_1-1)}{\Theta_1 l_1}}
\nonumber\\
&& \times
   \left(
      \int_{{\R^N}}
         \Phi_{1}(|u_{\lambda}|)|u_{\lambda}|^{m_1}
      dx
      +
      \int_{{\R^N}}
           (\Phi_{1}^{-1}(1))^{\Theta_1 l_1(\beta_1-1)-\frac{\Theta_1 (r_1-l_1)}{l_1-1}-m_1}|u_{\lambda}|^{m_1}
     dx
     \right)^{\frac{l_1-1}{\Theta_1 l_1}}
\nonumber\\
&&
   +\lambda M_7 \int_{{\R^N}} |u_{\lambda}|^{r_1-l_1}|w_k|^{l_1}dx\nonumber\\
& \leq &
    \lambda M_8 2^{r_2-1-\frac{(m_2+1)(l_1-1)(\Theta_1-1)}{\Theta_1 l_1}} 2^{(l_1-1)(\beta_1-1)-\frac{r_1-l_1}{l_1}-\frac{(m_1+1)(l_1-1)}{\Theta_1 l_1}}
   \left(\int_{{\R^N}} |w_k|^{l_1}|u_\lambda|^{r_1-l_1}dx\right)^{\frac{1}{l_1}}
\nonumber\\
&&  \times
    \bigg(
    \left(
       \int_{{\R^N}}
           \Phi_{2}^{\frac{\sigma_1}{\sigma_1-1}}(|v_{\lambda}|)
      dx
    \right)^{\frac{\sigma_1-1}{\sigma_1}}
    \left(
      \int_{{\R^N}}
            |v_{\lambda}|^{m_2\sigma_1}
      dx
    \right)^{\frac{1}{\sigma_1}}\nonumber\\
&&
   +
      \int_{{\R^N}}
(\Phi_{2}^{-1}(1))^{\frac{\Theta_1 l_1(r_2-1)}{(\Theta_1-1)(l_1-1)}-m_2} |v_{\lambda}|^{m_2}
      dx
     \bigg)^{\frac{(l_1-1)(\Theta_1-1)}{\Theta_1 l_1}}
\nonumber\\
&&  \times
   \bigg(
      \left(
       \int_{{\R^N}}
           \Phi_{1}^{\frac{\sigma_2}{\sigma_2-1}}(|u_{\lambda}|)
      dx
    \right)^{\frac{\sigma_2-1}{\sigma_2}}
    \left(
      \int_{{\R^N}}
            |u_{\lambda}|^{m_1\sigma_2}
      dx
    \right)^{\frac{1}{\sigma_2}}\nonumber\\
&&
     +
      \int_{{\R^N}}
(\Phi_{1}^{-1}(1))^{\Theta_1 l_1(\beta_1-1)-\frac{\Theta_1 (r_1-l_1)}{l_1-1}-m_1}|u_{\lambda}|^{m_1}
     dx
     \bigg)^{\frac{l_1-1}{\Theta_1 l_1}}
\nonumber\\
&&
   +\lambda M_7 \int_{{\R^N}} |u_{\lambda}|^{r_1-l_1}|w_k|^{l_1}dx\nonumber\\
& \leq &
   \lambda M_8 2^{r_2-1-\frac{(m_2+1)(l_1-1)(\Theta_1-1)}{\Theta_1 l_1}} 2^{(l_1-1)(\beta_1-1)-\frac{r_1-l_1}{l_1}-\frac{(m_1+1)(l_1-1)}{\Theta_1 l_1}}
   \left(\int_{{\R^N}} |w_k|^{l_1}|u_\lambda|^{r_1-l_1}dx\right)^{\frac{1}{l_1}}
\nonumber\\
&&  \times
    \bigg(
    \left(
       \int_{{\R^N}}
           \Phi_{2}^{\frac{\sigma_1}{\sigma_1-1}}(|v_{\lambda}|)
      dx
    \right)^{\frac{(\sigma_1-1)(l_1-1)(\Theta_1-1)}{\sigma_1\Theta_1 l_1}}
    \left(
      \int_{{\R^N}}
            |v_{\lambda}|^{m_2\sigma_1}
      dx
    \right)^{\frac{(l_1-1)(\Theta_1-1)}{\sigma_1\Theta_1 l_1}}
    \nonumber\\
&&
   +
     (\Phi_{2}^{-1}(1))^{r_2-1-\frac{m_2(l_1-1)(\Theta_1-1)}{\Theta_1 l_1}}
     \left(
      \int_{{\R^N}}
          |v_{\lambda}|^{m_2}
      dx
     \right)^{\frac{(l_1-1)(\Theta_1-1)}{\Theta_1 l_1}}
     \bigg)
\nonumber\\
&&  \times
   \bigg(
      \left(
       \int_{{\R^N}}
           \Phi_{1}^{\frac{\sigma_2}{\sigma_2-1}}(|u_{\lambda}|)
      dx
    \right)^{\frac{(\sigma_2-1)(l_1-1)}{\sigma_2\Theta_1 l_1}}
    \left(
      \int_{{\R^N}}
            |u_{\lambda}|^{m_1\sigma_2}
      dx
    \right)^{\frac{l_1-1}{\sigma_2\Theta_1 l_1}}
    \nonumber\\
&&
     +
     (\Phi_{1}^{-1}(1))^{(\beta_1-1)(l_1-1)-\frac{r_1-l_1}{l_1}-\frac{m_1(l_1-1)}{\Theta_1 l_1}}
     \left(
      \int_{{\R^N}}
        |u_{\lambda}|^{m_1}
     dx
     \right)^{\frac{l_1-1}{\Theta_1 l_1}}
     \bigg)
\nonumber\\
&&
   +\lambda M_7 \int_{{\R^N}} |u_{\lambda}|^{r_1-l_1}|w_k|^{l_1}dx\nonumber\\
& \leq &
    \lambda M_8 2^{r_2-1-\frac{(m_2+1)(l_1-1)(\Theta_1-1)}{\Theta_1 l_1}} 2^{(l_1-1)(\beta_1-1)-\frac{r_1-l_1}{l_1}-\frac{(m_1+1)(l_1-1)}{\Theta_1 l_1}}
   \left(\int_{{\R^N}} |w_k|^{l_1}|u_\lambda|^{r_1-l_1}dx\right)^{\frac{1}{l_1}}
\nonumber\\
&&      \times
       C_{0,2}^{\frac{m_2(l_1-1)(\Theta_1-1)}{\Theta_1 l_1}}\left\|v_{\lambda}\right\|_{\Phi_{2}}^{\frac{m_2(l_1-1)(\Theta_1-1)}{\Theta_1 l_1}}
     \bigg(
         (\Phi_{2}^{-1}(1))^{r_2-1-\frac{m_2(l_1-1)(\Theta_1-1)}{\Theta_1 l_1}}
         \nonumber\\
&&
         +
         \max\left\{   \left\|v_{\lambda}\right\|_{\Phi_{2}}^{\frac{l_2(\sigma_1-1)(l_1-1)(\Theta_1-1)}{\sigma_1\Theta_1 l_1} }, \left\|v_{\lambda}\right\|_{\Phi_{2}}^{\frac{m_2(\sigma_1-1)(l_1-1)(\Theta_1-1)}{\sigma_1\Theta_1 l_1}}\right\}
     \bigg)
\nonumber\\
&&   \times
    C_{0,1}^{\frac{m_1(l_1-1)}{\Theta_1 l_1}}\left\|u_{\lambda}\right\|_{\Phi_{1}}^{\frac{m_1(l_1-1)}{\Theta_1 l_1}}
    \bigg(
            (\Phi_{1}^{-1}(1))^{(\beta_1-1)(l_1-1)-\frac{r_1-l_1}{l_1}-\frac{m_1(l_1-1)}{\Theta_1 l_1}}
            \nonumber\\
&&
          +
            \max\left\{   \left\|u_{\lambda}\right\|_{\Phi_{1}}^{\frac{l_1(\sigma_2-1)(l_1-1)}{\sigma_2\Theta_1 l_1} }, \left\|u_{\lambda}\right\|_{\Phi_{1}}^{\frac{m_1(\sigma_2-1)(l_1-1)}{\sigma_2\Theta_1 l_1}}\right\}
     \bigg)
\nonumber\\
&&
   +\lambda M_7 \int_{{\R^N}} |u_{\lambda}|^{r_1-l_1}|w_k|^{l_1}dx\nonumber\\
& = &
    \lambda M_8 C_0
       C_{0,2}^{\frac{m_2(l_1-1)(\Theta_1-1)}{\Theta_1 l_1}}\left\|v_{\lambda}\right\|_{\Phi_{2}}^{\frac{m_2(l_1-1)(\Theta_1-1)}{\Theta_1 l_1}}
    C_{0,1}^{\frac{m_1(l_1-1)}{\Theta_1 l_1}}\left\|u_{\lambda}\right\|_{\Phi_{1}}^{\frac{m_1(l_1-1)}{\Theta_1 l_1}}
    \nonumber\\
&&    \times
     \left(
           \max\left\{   \left\|v_{\lambda}\right\|_{\Phi_{2}}^{\frac{l_2(\sigma_1-1)(l_1-1)(\Theta_1-1)}{\sigma_1\Theta_1 l_1} }, \left\|v_{\lambda}\right\|_{\Phi_{2}}^{\frac{m_2(\sigma_1-1)(l_1-1)(\Theta_1-1)}{\sigma_1\Theta_1 l_1}}\right\}
       +
         (\Phi_{2}^{-1}(1))^{r_2-1-\frac{m_2(l_1-1)(\Theta_1-1)}{\Theta_1 l_1}}
     \right)
\nonumber\\
&&      \times
    \left(
           \max\left\{   \left\|u_{\lambda}\right\|_{\Phi_{1}}^{\frac{l_1(\sigma_2-1)(l_1-1)}{\sigma_2\Theta_1 l_1} }, \left\|u_{\lambda}\right\|_{\Phi_{1}}^{\frac{m_1(\sigma_2-1)(l_1-1)}{\sigma_2\Theta_1 l_1}}\right\}
       +
         (\Phi_{1}^{-1}(1))^{(\beta_1-1)(l_1-1)-\frac{r_1-l_1}{l_1}-\frac{m_1(l_1-1)}{\Theta_1 l_1}}
     \right)
\nonumber\\
&&  \times
  \left(\int_{{\R^N}} |w_k|^{l_1}|u_\lambda|^{r_1-l_1}dx\right)^{\frac{1}{l_1}}
   +\lambda M_7 \int_{{\R^N}} |u_{\lambda}|^{r_1-l_1}|w_k|^{l_1}dx\nonumber\\
&\le  &   \lambda M_7 \int_{{\R^N}} |u_{\lambda}|^{r_1-l_1}|w_k|^{l_1}dx
        + \lambda M_8 C_1 \left(\int_{{\R^N}} |w_k|^{l_1}|u_\lambda|^{r_1-l_1}dx\right)^{\frac{1}{l_1}}
\end{eqnarray}
for some   $\sigma_1\in \left(1,\frac{l_2^*}{m_2}\right]$, $\sigma_2\in \left(1,\frac{l_1^*}{m_1}\right]$
 and any given $\beta_1\in I_\beta :=\left(1+\frac{r_1-l_1}{l_1(l_1-1)}+\frac{m_1+1}{\Theta_1 l_1},+\infty\right)$, where
$$
\Phi_{1}^{\Theta_1 l_1(\beta_1-1)-\frac{\Theta_1 (r_1-l_1)}{l_1-1}-m_1}(|u_{\lambda}|)<\Phi_{1}(|u_{\lambda}|)
$$
which holds by (\ref{2.2tj}).  Moreover, for the validity of the interval of $\Theta_1$, we need the following restriction
 \begin{eqnarray}\label{ccc1}
 r_2 > 1+\frac{(1+m_2)(l_1-1)(\Theta_1-1)}{l_1\Theta_1} \nonumber
\end{eqnarray}
to make
$$
\Phi_{2}^{\frac{\Theta_1 l_1(r_2-1)}{(\Theta_1-1)(l_1-1)}-m_2}(|v_{\lambda}|)<\Phi_{2}(|v_{\lambda}|)
$$
 hold by (\ref{2.2tj}).
In addition,   the range of $\sigma_1$ and $\sigma_2$ is determined by (\ref{2.5}), and we assume that
\begin{eqnarray}
C_1
&=& C_{0}
       C_{0,1}^{\frac{m_1(l_1-1)}{\Theta_1 l_1}}C_{0,2}^{\frac{m_2(l_1-1)(\Theta_1-1)}{\Theta_1 l_1}}
       \left\|v_{\lambda}\right\|_{\Phi_{2}}^{\frac{m_2(l_1-1)(\Theta_1-1)}{\Theta_1 l_1}}
      \left\|u_{\lambda}\right\|_{\Phi_{1}}^{\frac{m_1(l_1-1)}{\Theta_1 l_1}}
    \nonumber\\
&&
    \times \left(
           \max\left\{   \left\|v_{\lambda}\right\|_{\Phi_{2}}^{\frac{l_2(\sigma_1-1)(l_1-1)(\Theta_1-1)}{\sigma_1\Theta_1 l_1} }, \left\|v_{\lambda}\right\|_{\Phi_{2}}^{\frac{m_2(\sigma_1-1)(l_1-1)(\Theta_1-1)}{\sigma_1\Theta_1 l_1}}\right\}
       +
         (\Phi_{2}^{-1}(1))^{r_2-1-\frac{m_2(l_1-1)(\Theta_1-1)}{\Theta_1 l_1}}
     \right)
\nonumber\\
&&
   \times \left(
           \max\left\{   \left\|u_{\lambda}\right\|_{\Phi_{1}}^{\frac{l_1(\sigma_2-1)(l_1-1)}{\sigma_2\Theta_1 l_1} }, \left\|u_{\lambda}\right\|_{\Phi_{1}}^{\frac{m_1(\sigma_2-1)(l_1-1)}{\sigma_2\Theta_1 l_1}}\right\}
       +
         (\Phi_{1}^{-1}(1))^{(\beta_1-1)(l_1-1)-\frac{r_1-l_1}{l_1}-\frac{m_1(l_1-1)}{\Theta_1 l_1}}
     \right),\nonumber
\end{eqnarray}
where
$C_{0}=2^{r_2-1+(l_1-1)(\beta_1-1)-\frac{(m_2+1)(l_1-1)(\Theta_1-1)+(m_1+1)(l_1-1)}{\Theta_1 l_1}-\frac{r_1-l_1}{l_1}}.$
\par
By Gagliard-Nirenberg-Sobolev inequality, (\ref{kk2}) and H$\ddot{\mbox{o}}$lder inequality, we have
 \begin{eqnarray}\label{kkk2}
         \left(\int_{{\R^N}} |w_k|^{l_1^*}dx\right)^{\frac{l_1}{l_1^*}}
& \le & D_1 \int_{{\R^N}} |\nabla w_k|^{l_1}dx\nonumber\\
&\le  &\lambda D_2\beta_1^{l_1}\int_{{\R^N}} |u_{\lambda}|^{r_1-l_1}|w_k|^{l_1}dx+\lambda D_3 \beta_1^{l_1}  \left(\int_{{\R^N}} |u_{\lambda}|^{r_1-l_1}|w_k|^{l_1}dx\right)^{\frac{1}{l_1}}\nonumber\\
&\le  &\lambda D_2\beta_1^{l_1}\left(\int_{{\R^N}} |u_{\lambda}|^{l_1^*}\right)^{\frac{r_1-l_1}{l_1^*}}\left(\int_{{\R^N}}|w_k|^{\frac{l_1l_1^*}{l_1^*-r_1+l_1}}dx\right)^{\frac{l_1^*-r_1+l_1}{l_1^*}}\nonumber\\
&     &+\lambda D_3 \beta_1^{l_1}\left(\int_{{\R^N}} |u_{\lambda}|^{l_1^*}\right)^{\frac{r_1-l_1}{l_1l_1^*}}\left(\int_{{\R^N}}|w_k|^{\frac{l_1l_1^*}{l_1^*-r_1+l_1}}dx\right)^{\frac{l_1^*-r_1+l_1}{l_1l_1^*}}\nonumber\\
&\le  &\lambda D_2\beta_1^{l_1}\left(\int_{{\R^N}} |u_{\lambda}|^{l_1^*}\right)^{\frac{r_1-l_1}{l_1^*}}\left(\int_{{\R^N}}|u_\lambda|^{\frac{\beta_1l_1l_1^*}{l_1^*-r_1+l_1}}dx\right)^{\frac{l_1^*-r_1+l_1}{l_1^*}}\nonumber\\
&     &+\lambda D_3 \beta_1^{l_1}\left(\int_{{\R^N}} |u_{\lambda}|^{l_1^*}\right)^{\frac{r_1-l_1}{l_1l_1^*}}\left(\int_{{\R^N}}|u_\lambda|^{\frac{\beta_1l_1l_1^*}{l_1^*-r_1+l_1}}dx\right)^{\frac{l_1^*-r_1+l_1}{l_1l_1^*}}\nonumber\\
&\le  &\lambda D_2\beta_1^{l_1}\|u_\lambda\|_{1,\Phi_1}^{r_1-l_1}\left(\int_{{\R^N}}|u_\lambda|^{\frac{\beta_1l_1l_1^*}{l_1^*-r_1+l_1}}dx\right)^{\frac{l_1^*-r_1+l_1}{l_1^*}}\nonumber\\
&     &+\lambda D_3 \beta_1^{l_1}\|u_\lambda\|_{1,\Phi_1}^{\frac{r_1-l_1}{l_1}}\left(\int_{{\R^N}}|u_\lambda|^{\frac{\beta_1l_1l_1^*}{l_1^*-r_1+l_1}}dx\right)^{\frac{l_1^*-r_1+l_1}{l_1l_1^*}}.\nonumber
 \end{eqnarray}
Note that $r_1\in \bigg(\max\left\{m_1,m_2,1+\frac{(1+m_1)(l_2-1)(\Theta_2-1)}{l_2\Theta_2}\right\},
\min\left\{l_1^*,1+ \frac{l_1^*(l_1-1)}{l_1}-\frac{(l_1-1)(1+m_1)}{\Theta_1 l_1}\right\}\bigg)$. Let $\beta_1=1+\frac{l_1^*-r_1}{l_1}$. Then $\beta_1\in I_\beta$. Thus $\frac{\beta_1 l_1l_1^*}{l_1^*-r_1+l_1}=l_1^*$. Thus, we have
 \begin{eqnarray}\label{kkk3}
         \left(\int_{{\R^N}} |w_k|^{l_1^*}dx\right)^{\frac{l_1}{l_1^*}}
 &\le  &\lambda D_2\beta_1^{l_1}\|u_\lambda\|_{1,\Phi_1}^{r_1-l_1}\|u_\lambda\|_{\beta_1\alpha_1^*}^{l_1\beta_1}+\lambda D_3 \beta_1^{l_1}\|u_\lambda\|_{1,\Phi_1}^{\frac{r_1-l_1}{l_1}}\|u_\lambda\|_{\beta_1\alpha_1^*}^{\beta_1},
 \end{eqnarray}
 where $\alpha_1^*=\frac{l_1l_1^*}{l_1^*-r_1+l_1}$. Lemma 3.3 implies that there exists sufficient large $\Lambda_3>\Lambda_0$ such that
  $\|u_\lambda\|_{1,\Phi_1}<\min\{1,\frac{1}{C_{0,1}}\}$ for all $\lambda>\Lambda_3$ and then (\ref{2.5}) implies that
  $\|u_\lambda\|_{\beta_1\alpha_1^*}=\|u_\lambda\|_{l_1^*}<\min\{C_{0,1},1\}$. Since $r_1-l_1>\frac{r_1-l_1}{l_1}$ and $l_1\beta_1>\beta_1$. Hence, by (\ref{kkk3}), we have
   \begin{eqnarray}\label{kkk4}
         \left(\int_{{\R^N}} |w_k|^{l_1^*}dx\right)^{\frac{l_1}{l_1^*}}
 &\le  &\lambda(D_2+ D_3) \beta_1^{l_1}\|u_\lambda\|_{1,\Phi_1}^{\frac{r_1-l_1}{l_1}}\|u_\lambda\|_{\beta_1\alpha_1^*}^{\beta_1}.
 \end{eqnarray}
 Then it follows from the definition of $w_k$, Fatou's Lemma and (\ref{kkk4}) that
 \begin{eqnarray}\label{kkk5}
         \|u_\lambda\|_{\beta_1l_1^*}
 &\le  &\left(\lambda (D_2+ D_3) \beta_1^{l_1}\|u_\lambda\|_{1,\Phi_1}^{\frac{r_1-l_1}{l_1}}\right)^{\frac{1}{l_1\beta_1}}\|u_\lambda\|_{\beta_1\alpha_1^*}^{\frac{1}{l_1}}.\nonumber
 \end{eqnarray}
 Next, we start the iteration process. For each $n=0,1,2\cdots,$ we define $\beta_{1}^{(n+1)}=\frac{l_1^*}{\alpha_1^*}\beta_1^{(n)}$, where $\beta_1^{(0)}=\beta_1$.
  \begin{eqnarray}\label{kkk6}
         \|u_\lambda\|_{\beta_1^{(1)}l_1^*}
  &\le  &\left(\lambda (D_2+D_3) \beta_1^{(1)l_1}\|u_\lambda\|_{1,\Phi_1}^{\frac{r_1-l_1}{l_1}}\right)^{\frac{1}{l_1\beta_1^{(1)}}}\|u_\lambda\|_{\beta_1^{(1)}\alpha_1^*}^{\frac{1}{l_1}}\nonumber\\
 &\le  &\left(\lambda (D_2+D_3) \beta_1^{(1)l_1}\|u_\lambda\|_{1,\Phi_1}^{\frac{r_1-l_1}{l_1}}\right)^{\frac{1}{l_1\beta_1^{(1)}}}
 \left(\left(\lambda(D_2+D_3) \beta_1^{l_1}\|u_\lambda\|_{1,\Phi_1}^{\frac{r_1-l_1}{l_1}}\right)^{\frac{1}{l_1\beta_1}}
 \|u_\lambda\|_{\beta_1\alpha_1^*}^{\frac{1}{l_1}}\right)^{\frac{1}{l_1}}
 \nonumber\\
  &=  &(\lambda(D_2+D_3) \|u_\lambda\|_{1,\Phi_1}^{\frac{r_1-l_1}{l_1}})^{\frac{1}{l_1\beta_1^{(1)}}+\frac{1}{l_1^2\beta_1}}
 \beta_1^{\frac{1}{l_1\beta_1}}\beta_1^{(1)\frac{1}{\beta_1^{(1)}}}\|u_\lambda\|_{\beta_1\alpha_1^*}^{\frac{1}{l_1^2}}\nonumber
 \end{eqnarray}
 and then
 \begin{eqnarray}\label{kkk7}
         \|u_\lambda\|_{\beta_1^{(2)}l_1^*}
  &\le  &\left(\lambda (D_2+D_3) \beta_1^{(2)l_1}\|u_\lambda\|_{1,\Phi_1}^{\frac{r_1-l_1}{l_1}}\right)^{\frac{1}{l_1\beta_1^{(2)}}}\|u_\lambda\|_{\beta_1^{(2)}\alpha_1^*}^{\frac{1}{l_1}}\nonumber\\
 &\le  &
 \left((\lambda(D_2+D_3) \|u_\lambda\|_{1,\Phi_1}^{\frac{r_1-l_1}{l_1}})^{\frac{1}{l_1\beta_1^{(1)}}+\frac{1}{l_1^2\beta_1}}
 \beta_1^{\frac{1}{l_1\beta_1}}\beta_1^{(1)\frac{1}{\beta_1^{(1)}}}\|u_\lambda\|_{\beta_1\alpha_1^*}^{\frac{1}{l_1^2}}\right)^{\frac{1}{l_1}}\cdot\nonumber\\
&&
 \left(\lambda (D_2+D_3) \beta_1^{(2)l_1}\|u_\lambda\|_{1,\Phi_1}^{\frac{r_1-l_1}{l_1}}\right)^{\frac{1}{l_1\beta_1^{(2)}}}
 \nonumber\\
  &=  &\left(\lambda(D_2+D_3) \|u_\lambda\|_{1,\Phi_1}^{\frac{r_1-l_1}{l_1}}\right)^{\frac{1}{l_1\beta_1^{(2)}}+\frac{1}{l_1^2\beta_1^{(1)}}+\frac{1}{l_1^3\beta_1}}
 \beta_1^{\frac{1}{l_1^2\beta_1}}\beta_1^{(1)\frac{1}{l_1\beta_1^{(1)}}}\beta_1^{(2)\frac{1}{\beta_1^{(2)}}}\|u_\lambda\|_{\beta_1\alpha_1^*}^{\frac{1}{l_1^3}}\nonumber.
 \end{eqnarray}
 Repeating such process, we obtain that
 \begin{eqnarray}\label{kkk8}
   &    &      \|u_\lambda\|_{\beta_1^{(n)}l_1^*}\nonumber\\
  &\le  & \left(\lambda(D_2+D_3) \|u_\lambda\|_{1,\Phi_1}^{\frac{r_1-l_1}{l_1}}\right)^{\sum_{i=0}^n\frac{1}{l_1^{n+1-i}\beta_1^{(i)}}}
 \Pi_{i=0}^n\beta_1^{(i)\frac{1}{l_1^{n-i}\beta_1^{(i)}}}\|u_\lambda\|_{\beta_1\alpha_1^*}^{\frac{1}{l_1^{n+1}}}\nonumber\\
 &= & D_4\left(\lambda(D_2+D_3) \|u_\lambda\|_{1,\Phi_1}^{\frac{r_1-l_1}{l_1}}\right)^{\sum_{i=0}^n\frac{1}{\beta_1l_1l_1^{n-i}\left(\frac{l_1^*}{\alpha_1^*}\right)^{i}}}
\|u_\lambda\|_{\beta_1\alpha_1^*}^{\frac{1}{l_1^{n+1}}},
 \end{eqnarray}
where
$D_4= \beta_1^{\frac{1}{\beta_1}\sum_{i=0}^n\left(\frac{1}{l_1}\right)^{n-i}\left(\frac{\alpha_1^*}{l_1^*}\right)^{i}}
        \left(\frac{l_1^*}{\alpha_1^*}\right)^{\frac{1}{\beta_1}\sum_{i=0}^n i\left(\frac{1}{l_1}\right)^{n-i}\left(\frac{\alpha_1^*}{l_1^*}\right)^{i}}.$
\par
 Notice that
 $$
  \alpha_1^*<l_1^*,\;\;\;  \lim_{n\to \infty}\sum_{i=0}^n\frac{1}{\beta_1l_1l_1^{n-i}\left(\frac{l_1^*}{\alpha_1^*}\right)^{i}}=0,\ \  \
   \lim_{n\to \infty}\sum_{i=0}^n i\left(\frac{1}{l_1}\right)^{n-i}\left(\frac{\alpha_1^*}{l_1^*}\right)^{i}=0,\ \  \
  \lim_{n\to \infty}\frac{1}{l_1^{n+1}}=0.
 $$
 So (\ref{kkk8}) implies that
 $
         \|u_\lambda\|_{\infty}\le 1
$
 for all $\lambda>\Lambda_3$.
\par

Similarly, there exists $\Lambda_4>\Lambda_0$ such that $  \|v_\lambda\|_{\infty}\le 1$ for all  $\lambda>\Lambda_4$. Let $\Lambda_*=\max\{\Lambda_3,\Lambda_4\}$. Then we complete the proof.
\qed
 \vskip2mm
 \noindent
{\bf Proof of Theorem 1.1.}
 By  Lemma 3.4, for each $\lambda>\Lambda_*$, we have
 $$ \|(u_{\lambda},v_{\lambda})\|_{\infty}\le \|u_{\lambda}\|_{\infty}+\|v_{\lambda}\|_{\infty} \le 2$$ which implies that
$\widetilde{F}(x,u_{\lambda},v_{\lambda})={F}(x,u_{\lambda},v_{\lambda})$ for all $x\in \R^N$. Hence, $(u_{\lambda},v_\lambda)$ is a nontrivial weak solution of system (\ref{aa1}) and Lemma 3.3 implies that  $\|u_\lambda\|_{1,\Phi_1}\to 0$ and $\|v_\lambda\|_{1,\Phi_2}\to 0$ as $\lambda\to\infty$.
 \qed

 \vskip2mm
 {\section{Example}}
   \setcounter{equation}{0}
\begin{eqnarray}\label{d1}
 \begin{cases}
  -\mbox{div}[(4|\nabla u|^{2}+5|\nabla u|^{3})\nabla u]+\left(1+\sum_{i=1}^{6}\cos^2\pi x_i\right)(4|u|^{2}+5|u|^{3})u\\
 \;\;\; =\lambda F_u(x,u,v), \ \ x\in \R^6,\\
   -\mbox{div}\left[\left(4|\nabla v|^{2}\log(2+|\nabla v|)+\dfrac{|\nabla v|^{3}}{1+|\nabla v|}\right)\nabla v\right]+\left(1+\sum_{i=1}^{6}\sin^2\pi x_i\right)\left(4|v|^{2}\log(1+| v|)+\dfrac{|v|^{3}}{1+|v|}\right)v\\
  \;\;\;=\lambda F_v(x,u,v), \ \ x\in \R^6,\\
 \end{cases}
\end{eqnarray}
where
\begin{eqnarray}\label{dd3}
F(x,t,s)= \sigma(t,s)b(x)\left(|t|^{\frac{17}{2}}+|s|^{\frac{17}{2}}+|t|^{7}|s|^{7}\right)+(1-\sigma(t,s))b(x)\left(|t|^3+|s|^3\right)
\end{eqnarray}
 for all $(x,t,s)\in \R^{N}\times \R\times \R$  with $\sigma(t,s)$ defined by (\ref{dd2}) and
  $b(x)= \left(1+\sum_{i=1}^{6}\cos^2\pi x_i\right)$ (or $b(x)\equiv 1$).
 So,
\begin{eqnarray*}
 F(x,t,s)= \begin{cases}
           b(x)\left(|t|^{\frac{17}{2}}+|s|^{\frac{17}{2}}+|t|^{7}|s|^{7}\right), \;\;\; \;\;\; \;\;\;\;\;\;\; \;\;\;\;\;\;\; \;\; \;\;\;\;\; \;\; \;\;\; \;
                                                                                              \text { if }\;\;|(t,s)| \le 4, \\
             \sin\dfrac{\pi(t^2+s^2-64)^2}{4608}b(x)\left(|t|^{\frac{17}{2}}+|s|^{\frac{17}{2}}+|t|^{7}|s|^{7}\right)\\
             +\left(1-\sin\dfrac{\pi(t^2+s^2-64)^2}{4608}\right)b(x)(|t|^{3}+|s|^{3}), \;\;\; \;\;\; \;\; \;             \text { if }\;\;4 <|(t,s)| \leq 8, \\
           b(x)(|t|^{3}+|s|^{3}), \;\;\; \;\;\; \;\;\; \;\;\; \;\;\; \;\;\; \;\;\; \;\;\; \;\;\;\;\;\;\; \;\;\;\; \;\;\;\;\; \;\; \;\;\; \;\;\; \;\;\;\; \;\;
                                                                                              \text { if }\;\;|(t,s)| > 8,
          \end{cases}
   \end{eqnarray*}
   and we also draw the figure of $F$ (see figure (e)-(h) below). 
\begin{figure}[htbp]
\centering
\subfigure[$(t,s)\in((-0.5,0.5),(-0.5,0.5))$]{
\includegraphics[width=5.5cm]{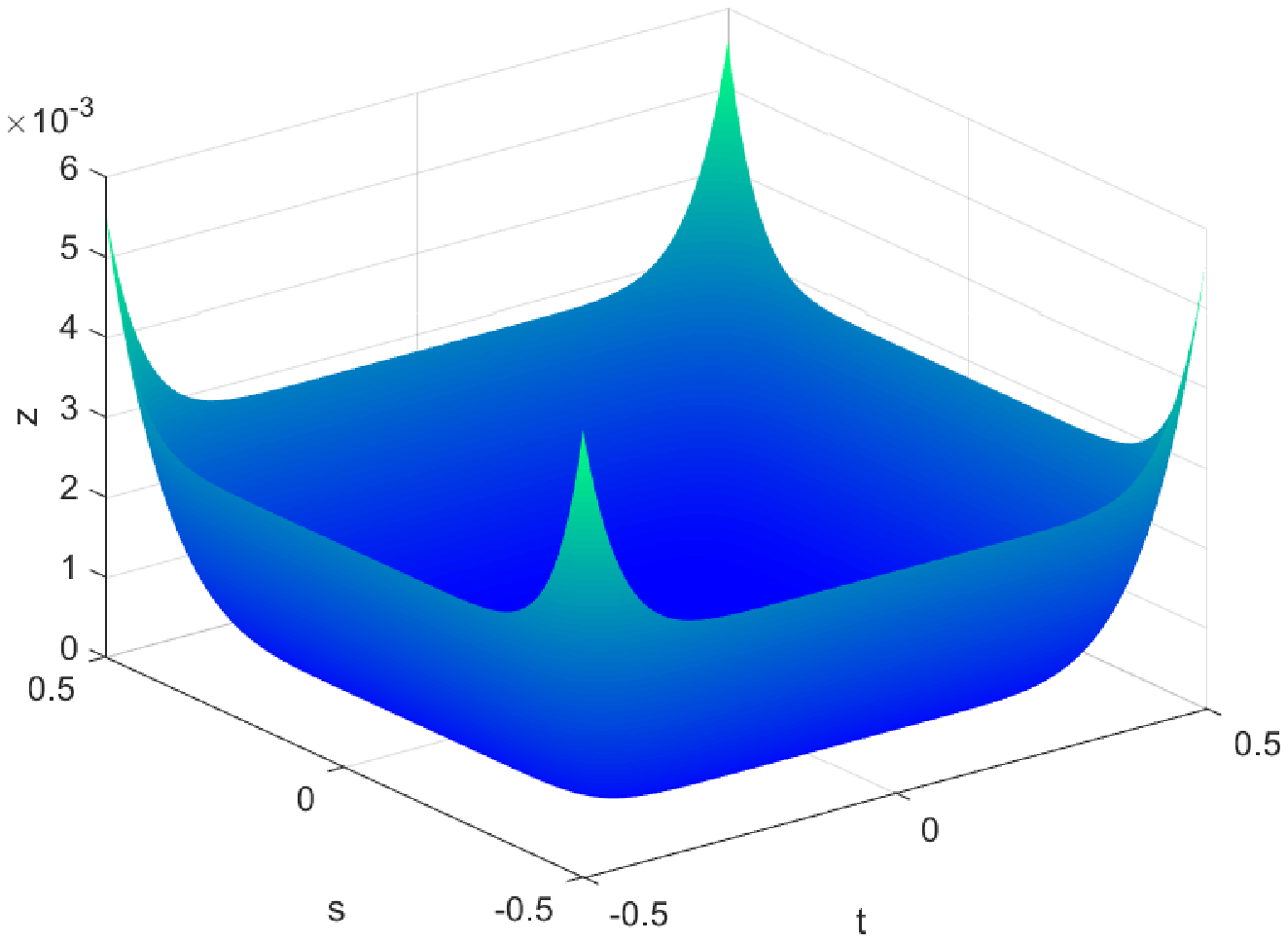}
}
\quad
\subfigure[$(t,s)\in((-4,4),(-4,4))$]{
\includegraphics[width=5.5cm]{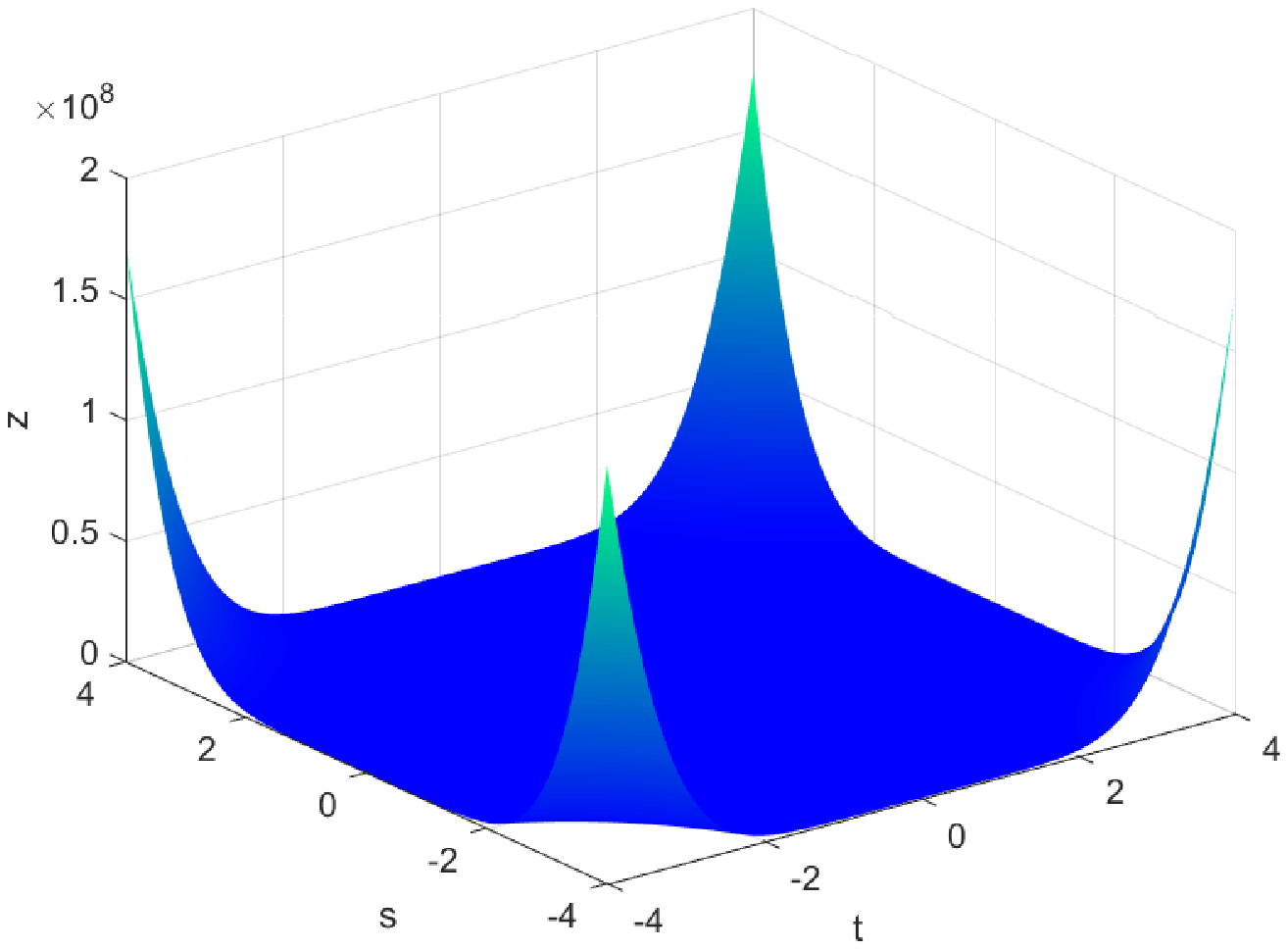}
}
\quad
\subfigure[$(t,s)\in((-8,8),(-8,8))$]{
\includegraphics[width=5.5cm]{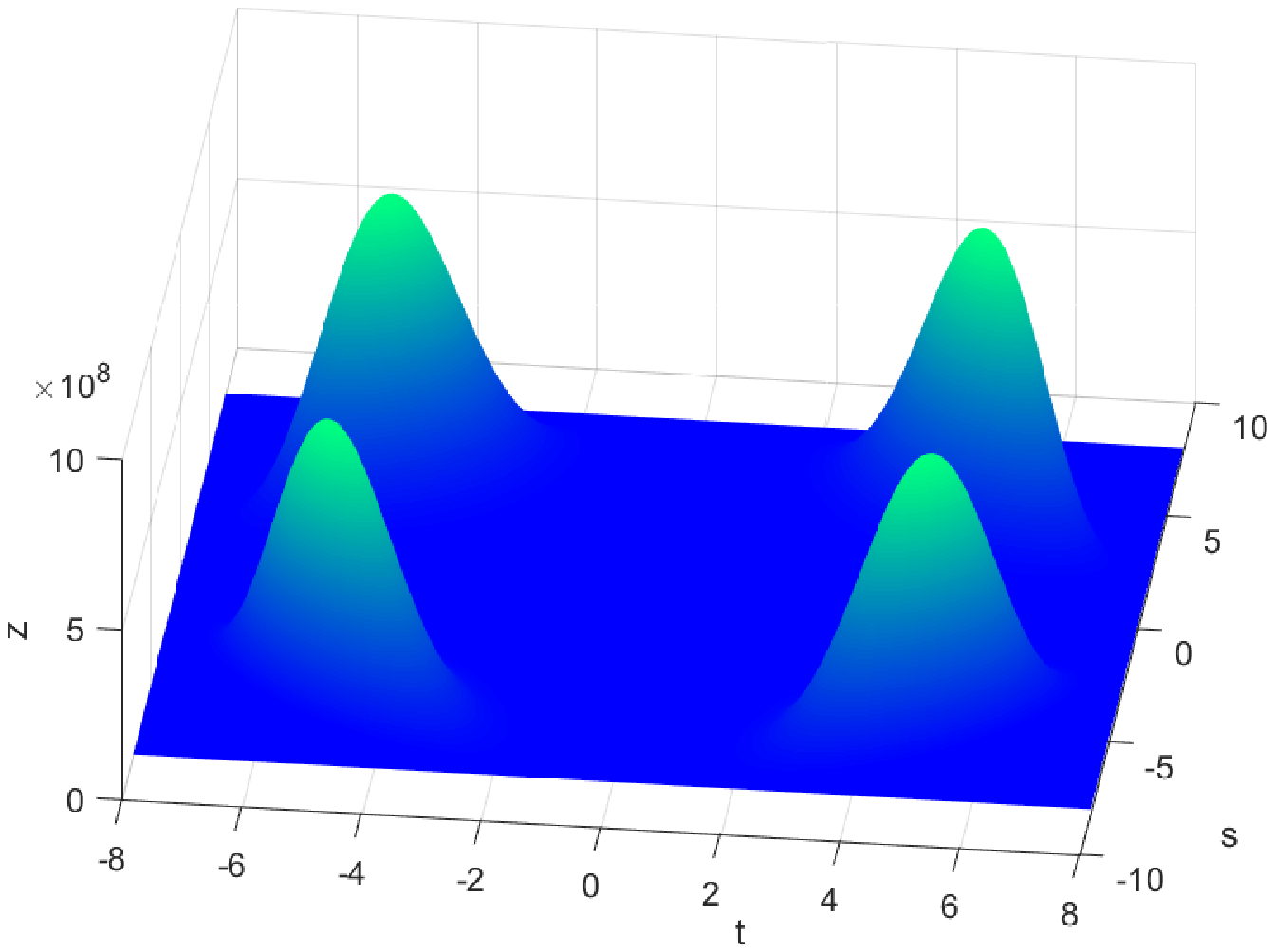}
}
\quad
\subfigure[$(t,s)\in((-2000,2000),(-2000,2000))$]{
\includegraphics[width=5.5cm]{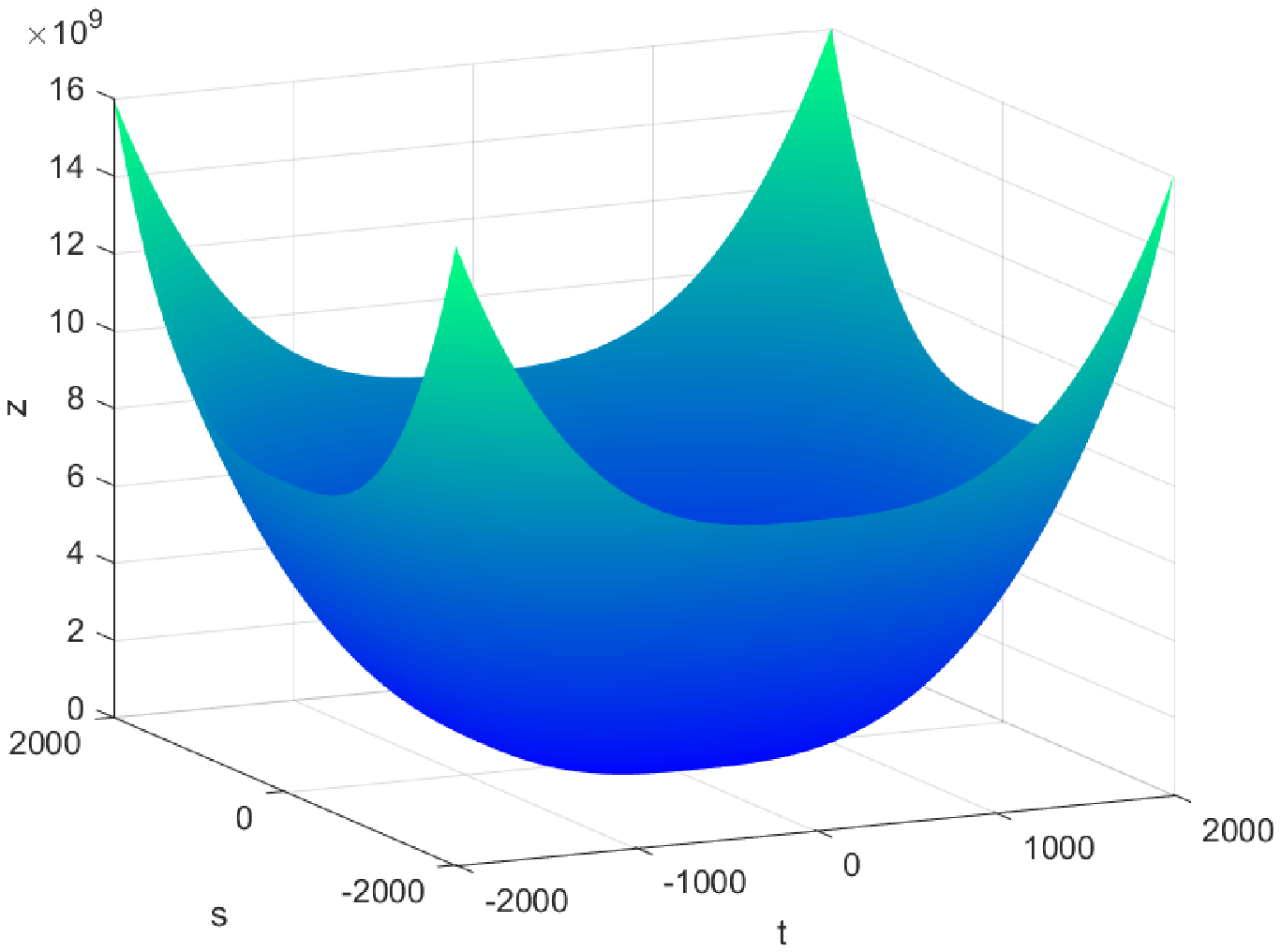}
}

\end{figure}
   \par
Let $N=6$, $\phi_1(t)=4|t|^{2}+5|t|^{3}$ and $\phi_2(t)=4|t|^{2}\log(2+| t|)+\dfrac{|t|^{3}}{1+|t|}$. Then $\phi_i (i=1,2)$ satisfy $(\phi_1)$-$(\phi_3)$ and $l_{1}=l_{2}=4$, $m_{1}=m_{2}=5$ (see \cite{Wang2016}). So, $l^{*}_{1}=l^{*}_{2}=12$.
\par
Let $V_1(x)=1+\sum_{i=1}^{6}\cos^2\pi x_i$ and $V_2(x)=1+\sum_{i=1}^{6}\sin^2\pi x_i$
 for all $(x,t,s)\in \R^{N}\times \R\times \R$. Then it is obvious that $V_i$ satisfies (V0) and (V1).
 \par
  Let
 \begin{eqnarray*}
 A(t,s)& :=  & \dfrac{t\pi(t^2+s^2-64)}{1152}\cos\dfrac{\pi(t^2+s^2-64)^2}{4608}b(x)\left(|t|^{\frac{17}{2}}+|s|^{\frac{17}{2}}+|t|^{7}|s|^{7}\right)\\
  &  &  +\sin\dfrac{\pi(t^2+s^2-64)^2}{4608}b(x)\left(\dfrac{17}{2}|t|^{\frac{13}{2}}t+7|t|^{5}|s|^{7}t\right)\\
  &   &  -\dfrac{t\pi(t^2+s^2-64)}{1152}\cos\dfrac{\pi(t^2+s^2-64)^2}{4608}b(x)\left(|t|^3+|s|^3\right)
         +3\left(1-\sin\dfrac{\pi(t^2+s^2-64)^2}{4608}\right)a(x)|t|t.
 \end{eqnarray*}
 and
  \begin{eqnarray*}
 B(t,s)& :=  & \dfrac{s\pi(t^2+s^2-64)}{1152}\cos\dfrac{\pi(t^2+s^2-64)^2}{4608}b(x)\left(|t|^{\frac{17}{2}}+|s|^{\frac{17}{2}}+|t|^{7}|s|^{7}\right)\\
  &  &  +\sin\dfrac{\pi(t^2+s^2-64)^2}{4608}b(x)\left(\dfrac{17}{2}|s|^{\frac{13}{2}}s+7|s|^{5}|t|^{7}s\right)\\
  &   &  -\dfrac{s\pi(t^2+s^2-64)}{1152}\cos\dfrac{\pi(t^2+s^2-64)^2}{4608}b(x)\left(|t|^3+|s|^3\right)
         +3\left(1-\sin\dfrac{\pi(t^2+s^2-64)^2}{4608}\right)b(x)|s|s.
 \end{eqnarray*}
 It is easy to see that $F$ satisfies (F0) and
 \begin{eqnarray*}
          F_t(x,t,s)
  &  =  & \begin{cases}
   b(x)\left(\dfrac{17}{2}|t|^{\frac{13}{2}}t+7|t|^{5}|s|^{7}t\right),
    & \mbox{ if }  |(t,s)|\le 4 ,\\
        A(t,s), & \text { if }\;\; 4<|(t,s)|\le 8\\
    3b(x)|t|t,
           & \mbox{ if }  |(t,s)|> 8,\\
    \end{cases}\\
         F_s(x,t,s)
  &  =  &\begin{cases}
  b(x)\left(\dfrac{17}{2}|s|^{\frac{13}{2}}s+7|s|^{5}|t|^{7}s\right), & \mbox{ if } |(t,s)|\le 4, \\
    B(t,s), & \text { if }\;\; 4< |(t,s)|\le 8\\
 3b(x)|s|s,
  & \mbox{ if }
          |(t,s)|>8.
  \end{cases}
\end{eqnarray*}
Thus, we  get
\begin{eqnarray*}
 F(x,t,s)&=&b(x)\left(|t|^{\frac{17}{2}}+|s|^{\frac{17}{2}}+|t|^{7}|s|^{7}\right)\\
&  \ge & \begin{cases}
           |t|^{9}+|s|^{9}, &  \text { if }\;\;|(t,s)| < 1, \\
        \dfrac{1}{16}(|t|^2+|s|^2)\left(|t|^{7}+|s|^{7}\right),&  \text { if }\;\; 1 \le |(t,s)|\le 4, \\
          \end{cases}\\
&\geq &  \dfrac{1}{16}|t|^{9}+\dfrac{1}{16}|s|^{9}, \;\;\;\mbox{ for all } x\in \R^{N} \mbox{ and } |(t,s)|\le 4.
\end{eqnarray*}
Then (F1) holds with $k_1=k_2=9\in (m_1,l_{1}^{*})=(m_2,l_{2}^{*})=\left(5,12\right)$, $M_1=M_2=\dfrac{1}{16}$. We also have
 \begin{eqnarray*}
|F_t(x,t,s)|&=& b(x)\left(\dfrac{17}{2}|t|^{\frac{15}{2}}+7|t|^{6}|s|^{7}\right)\\
& \le &  \begin{cases}
     \left(\dfrac{17}{2}+7\right)\times 7(|t|^{6}+|s|^{6}),  & \mbox{ if }  |(t,s)|<1\\
     \left(\dfrac{17}{2}+7\right)\times 7(t^{2}+s^{2})^4(|t|^{6}+|s|^{6}), & \mbox{ if }  1\le |(t,s)|\le 4
    \end{cases}\\
&\leq &   2^{41}(|t|^{6}+|s|^{6}), \;\;\;\mbox{ for all } x\in \R^{N} \mbox{ and } |(t,s)|\le 4.
\end{eqnarray*}
Similarly,
 \begin{eqnarray*}
|F_s(x,t,s)|&\leq &  2^{41}(|t|^{6}+|s|^{6}), \;\;\;\mbox{ for all } x\in \R^{N} \mbox{ and } |(t,s)|\le 4.
\end{eqnarray*}
Then (F2) holds with $M_3=M_4=2^{41}$ and
$r_1=7 \in \bigg(\max\left\{m_1,m_2,1+\frac{(1+m_1)(l_2-1)(\Theta_2-1)}{l_2\Theta_2}\right\},
1+ \frac{l_1^*(l_1-1)}{l_1}-\frac{(l_1-1)(1+m_1)}{\Theta_1 l_1}\bigg)=\bigg(\max\left\{5,1+\frac{9(\Theta_2-1)}{2\Theta_2}\right\},
10-\frac{9}{2\Theta_1}\bigg)$ and
$r_2=7 \in \bigg(\max\left\{m_1,m_2,1+\frac{(1+m_2)(l_1-1)(\Theta_1-1)}{l_1\Theta_1}\right\},
1+ \frac{l_2^*(l_2-1)}{l_2}-\frac{(l_2-1)(1+m_2)}{\Theta_2 l_2}\bigg)=\bigg(\max\left\{5,1+\frac{9(\Theta_1-1)}{2\Theta_1}\right\},
10-\frac{9}{2\Theta_2}\bigg)$  for some  $\Theta_1, \Theta_2>1$. In particular, taking  $\Theta_1=\Theta_2=6$, it is easy to see that
$r_1=r_2=7 \in \bigg(5, \dfrac{37}{4}\bigg)$.
Note that
\begin{eqnarray*}
 b(x)\left(|t|^{\frac{17}{2}}+|s|^{\frac{17}{2}}+|t|^{7}|s|^{7}\right)
\leq
b(x)\left(|t|^{\frac{17}{2}}+|s|^{\frac{17}{2}}+\frac{28}{17}|t|^{7}|s|^{7}\right).
\end{eqnarray*}
So
$$
F(x,t,s)\leq \frac{2}{17}tF_t(x,t,s)+\frac{2}{17}sF_s(x,t,s)=\frac{1}{\mu_1}tF_t(x,t,s)+\frac{1}{\mu_2}sF_s(x,t,s)
$$
for all $x\in \R^N$ and $|(t,s)|\le 4$, where $\mu_1=\mu_2=\dfrac{17}{2}>5=m_1=m_2$.
Thus, we have verified that system (\ref{d1}) satisfies all the conditions of Theorem 1.1. Hence, there exists $\Lambda_*>0$ such that system (\ref{aa1})  has  a nontrivial solution $(u_\lambda,v_\lambda)$ with $\|(u_\lambda,v_\lambda)\|_\infty\le 2$ for each $\lambda> \Lambda_*$ and $\|(u_\lambda,v_\lambda)\|\to 0$ as $\lambda\to\infty$.

\par
Next, we show that (\ref{dd3}) does not satisfy  Theorem A.
In fact, for (\ref{dd3}), we have
\begin{eqnarray*}
 &    &  \frac{1}{\mu_1}tF_t(x,t,s)+\frac{1}{\mu_2}sF_s(x,t,s)\\
 &\le &  \frac{1}{5}\left[tF_t(x,t,s)+sF_s(x,t,s)\right]\\
 &\le &  \frac{3}{5} b(x)|t|^3+\frac{3}{5} b(x)|s|^3\\
 & < &  b(x)|t|^3+b(x)|s|^3
\end{eqnarray*}
for all $x\in \R^{N}$, $|(t,s)|>4$ and any $\mu_i\in (5,+\infty)$, $i=1,2$, which shows that $F$ dose not satisfy (H2) in Theorem A.

 \vskip2mm
 {\section{A result for elliptic equation}}
   \setcounter{equation}{0}
 Notice that (F1)-(F3) hold for all $|(t,s)|\le 2$, which are stronger than $(f1)$-$(f3)$ where $|t|\le \delta$ for some positive constant $\delta$. This is because we can not obtain an estimate for $\|u\|_\infty$ and $\|v\|_\infty$ like Lemma 2.6 in  \cite{Medeiros} where $\|u\|_\infty$ grows with $\|u\|$. Instead, we obtain $\|u\|_\infty\le 1$ and $\|v\|_\infty\le 1$   (see Lemma 3.4), which are caused by the relation of $u$ and $v$ from the proof of Lemma 3.4.
 Hence, in this section, for the following quasilinear scalar equation
\begin{eqnarray}\label{741}
 \begin{cases}
  -\mbox{div}(\phi(|\nabla u|)\nabla u)+V(x)\phi(|u|)u=\lambda f(x, u), \ \ x\in \R^N,\\
  u\in W^{1,\Phi}(\R^N),
   \end{cases}
 \end{eqnarray}
where $N\ge 2$ and $\lambda>0$, we shall obtain a result similar to Theorem 1 in \cite{costa} and Theorem 1.2 in \cite{Medeiros}. To be precise, we have the following theorem:
\vskip2mm
\noindent
{\bf Theorem 5.1.} {\it Assume that $\phi$ satisfies $(\phi_1)$,  $(\phi_2)$ and the following conditions hold: \\
 $(\phi_3)$ there exists a positive constant $q$ such that $t^2\phi(|t|)\ge q|t|^{l}$ for all  $t\in \R$;\\
 (V0)\ $V\in C(\R^N,\R^+)$, $V_{\infty}:=\inf_{\R^N}V(x)>0$ and $V$ is a $1-$periodic function;\\
 (C1)  $f\in C(\R^N\times\R,\R)$ and $f$ is $1-$periodic in $x\in \R^N$;\\
(C2)\ there exist $\delta>0$, $k\in (m,K)$ and  $D_1>0$ such that
$$
F(x, t)\ge D_1 |t|^{k}
$$
for all $t\in (-\delta,\delta)$  and $x\in  \R^N$, where
$
K=\min\left\{l^{*},\frac{ml-l+l^*}{m}\right\}
$ and $F(x,t)=\int_{0}^{t}f(x,\xi)d\xi$;\\
(C3)\ there exist  $r\in (m, l^*)$ and  $D_2>0$  such that
\begin{eqnarray*}
   |f(x,t)|\le D_2 |t|^{r-1}
\end{eqnarray*}
for all $t\in (-\delta,\delta)$ and $x\in  \R^N$;\\
(C4) there exists $\mu>m$ such that
$$
0< F(x,t)\le \frac{1}{\mu}tf(x,t)
$$
for all $x\in  \R^N$ and  $t\in (-\delta,\delta)$ with $t\not=0$.\\
Then there exists $\Lambda_0>0$ such that equation (\ref{741})  has at least nontrivial solution $u_\lambda$ for each $\lambda> \Lambda_0$, $\|u_\lambda\|\to 0$  and $\|u_\lambda\|_\infty\to 0$ as $\lambda\to\infty$. }

 \vskip2mm
 The proof is similar to Theorem 1.1. The main difference is the result of Moser iteration. Next, we  outline the proof. We work on the space $W^{1,\Phi}(\R^N)$ with the norm $\|\cdot\|_{1,\Phi}$. Define the cut-off function $\rho\in C^1(\R,[0,1])$
  \begin{eqnarray*}
 \rho(t)= \begin{cases}
           1, \;\;\;  \text { if }\;\;|t| \leqslant \delta/2, \\
           0, \;\;\;  \text { if }\;\;|t| \geqslant  \delta
          \end{cases}
   \end{eqnarray*}
for all $t\in \R$  and  $t\rho'(t)\le 0 $. Similar to  (\ref{d3})-(\ref{d6}), some examples of $ \rho(t)$ can  also be given, for example,
 \begin{eqnarray*}
 \rho(t)= \begin{cases}
           1, &  \text { if }\;\;|t| < \frac{\delta}{2}, \\
          \sin\dfrac{8\pi(t^2-\delta^2)^2}{9\delta^4},    &  \text { if }\;\; \frac{\delta}{2}\le |t|\le \delta, \\
           0,&  \text { if }\;\;|t|>\delta.
          \end{cases}
 \end{eqnarray*}

 \vskip2mm
Let
$$
\widetilde{F}(x,t)=\rho(t)F(x,t)+(1-\rho(t)) D_3|t|^{r}.
$$
 \vskip2mm
 \noindent
   {\bf Lemma 5.1.} {\it Assume that (C1)-(C4) hold. Then \\
    (C1)$'$ $\widetilde{F}\in C^1(\R^N\times \R,\R)$, $\widetilde{F}$ is $1-$periodic in $x\in \R^N$ and $\widetilde{F}(x,0)=0$ for all $x\in \R^N$;\\
  (C2)$'$
  \begin{eqnarray*}
 0\le\widetilde{F}(x,t)\leq D_3 |t|^{r},  \ \ \mbox{for all }t\in \mathbb {R} \mbox{ and }x\in \R^N;
  \end{eqnarray*}
 (C3)$'$ there exists $D_4>0$  such that
 \begin{eqnarray*}
 |\widetilde{f}(x,t)|\le D_4|t|^{r-1}
   \end{eqnarray*}
 for all $t\in \mathbb {R}$ and $x\in \R^N$;\\
 (C3)$'$
   \begin{eqnarray*}
 \theta \widetilde{F}(x,t)\leq \widetilde{f}(x,t)t,  \ \ \mbox{for all }t\in \R/\{0\} \mbox{ and }x\in \R^N,
  \end{eqnarray*}
  where $\theta=\min\{r ,\mu \}$.
 }

 \vskip2mm
  Consider the modified problem
   \begin{eqnarray}\label{741-1}
\begin{cases}
  -\mbox{div}(\phi(|\nabla u|)\nabla u)+V(x)\phi(|u|)u=\lambda  \widetilde{f}(x, u), \ \ x\in \R^N,\\
  u\in W^{1,\Phi}(\R^N).
   \end{cases}
 \end{eqnarray}
 Define the functional $\widetilde{J}_{\lambda}: W^{1,\Phi}(\R^N)\to \R$ by
  $$
  \widetilde{J}_{\lambda}(u)=\int_{\R^N}\Phi(|\nabla u|)dx+\int_{\R^N}V(x)\Phi(|u|)dx-\lambda\int_{\R^N} \widetilde{F}(x,u)dx.
  $$
 It is easy to see that $ \widetilde{J}_{\lambda} $ is well defined and $\widetilde{J}_{\lambda}\in C^1(W^{1,\Phi}(\R^N),\R)$ and
\begin{eqnarray*}
          \langle \widetilde{J}'_{\lambda}(u),\tilde{u}\rangle
   =   \int_{\R^N}(\phi(|\nabla u|)\nabla u,\nabla \tilde{u})dx+\int_{\R^N}V(x)\phi(|u|)u\tilde{u}dx
       -\lambda\int_{\R^N}\widetilde{F}_u(x,u)\tilde{u}dx.
 \end{eqnarray*}
 \vskip2mm
 \noindent
   {\bf Lemma 5.2.} {\it $\widetilde{J}_{\lambda}$ satisfies the mountain pass geometry, that is,\\
   (i) there exist two positive constants $\gamma,\eta$ such that $\widetilde{J}_{\lambda}(u)\ge \eta$ for all $\|u\|=\gamma$; \\
   (ii) there exists $u_0\in C_0^\infty(\R^N)/\{0\}$ with $u_0>0$ and  $0<\|u_0\|_{\infty}<\frac{\delta}{2}$ such that $\widetilde{J}_{\lambda}(u_0)<0$.}
   \vskip2mm
   \noindent
{\bf Proof.}  Notice that (C1)$'$-(C2)$'$ imply those conditions of Theorem 1.5 in  \cite{Alves2014}. Then the proof is completed easily. \qed
 \vskip2mm
\par
 By Lemma 4.1-Lemma 4.3 in \cite{Alves2014} that system (\ref{741}) has a nontrivial solution $u_{\lambda}$ such that $\widetilde{J}_{\lambda}(u_{\lambda})=c_\lambda$ with
 \begin{eqnarray*}\label{kkk1}
 c_{{\lambda}}:=\inf_{\gamma\in\Gamma}\max_{t\in[0,1]}\widetilde{J}_{\lambda}(\gamma(t)),
\end{eqnarray*}
and
 $$
 \Gamma:=\{\gamma\in([0,1],W^{1,\Phi}(\R^N)):\gamma(0)=(0),\gamma(1)=u_0\}.
 $$

\vskip2mm
 \noindent
{\bf Lemma 5.3.} For each $\lambda>0$,  there exists $C_{**}>0$ such that
\begin{eqnarray*}
   \|u_{\lambda}\|_{1,{\Phi}}\leq  C_{**}\max\left\{\lambda^{-\frac{1}{k-l}},\lambda^{-\frac{l}{m(k-l)}}
\right\}.
 \end{eqnarray*}
 {\bf Proof.} It is easy to obtain the conclusion from the proof of Lemma 3.3.
 \vskip2mm
 \noindent
 {\bf Lemma 5.4.} {\it For each $\lambda>0$,   there exist  positive constants $C$ which only depends on $q, r,l,N$, respectively, such that}
\begin{eqnarray*}
&  &\|u_{\lambda}\|_\infty\le C(\lambda\|u_{\lambda}\|_{1,\Phi}^{r-l})^{\frac{1}{l^*-r}}\|u_{\lambda}\|_{1,\Phi}.
\end{eqnarray*}
{\bf Proof.} The proof is a direct generalization of  Lemma 2.6 in \cite{Medeiros}.    Since $u_{\lambda}$ is a critical point of $\widetilde{J}_{\lambda}$, we have
\begin{eqnarray}\label{752}
  \int_{{\R^N}}(\phi(|\nabla u_{\lambda}|)\nabla u_{\lambda},\nabla \tilde{u})dx+ \int_{{\R^N}}V(x)\phi(|u_{\lambda}|)u_{\lambda}\tilde{u}dx=\lambda \int_{{\R^N}} \widetilde{f}(x,u_{\lambda})\tilde{u}dx
\end{eqnarray}
for all $\tilde{u}\in W^{1,\Phi}(\R^N)$.  Without loss of generality, for each $k>0$, define
$$
u_k=\begin{cases}
  u_{\lambda},  & \mbox{if } u_{\lambda} \le k,\\
  k,  & \mbox{if } u_{\lambda}  > k,
\end{cases}
$$
$\varphi_k=|u_k|^{l(\beta-1)}u_\lambda$ and $w_k=u_{\lambda}|u_k|^{\beta-1}$ with $\beta>1$.
By $(\phi_4)$, we have $\phi(|\nabla u_{\lambda}|)|\nabla u_{\lambda}|^2\ge q |\nabla u_{\lambda}|^{l}$. Then taking $\tilde{u}=\varphi_k$ in (\ref{752}), by the definition of $u_k$ and (C3)$'$, we have
 \begin{eqnarray}
&     &      q \int_{{\R^N}}|\nabla u_{\lambda}|^{l} |u_k|^{l(\beta-1)}dx\nonumber\\
& \le &  \int_{{\R^N}}\phi(|\nabla u_{\lambda}|)\nabla u_{\lambda},\nabla \varphi_k)dx -l(\beta-1) \int_{{\R^N}}|u_k|^{l(\beta-1)-2}u_k u_\lambda\phi(|\nabla u_{\lambda}|)(\nabla u_{\lambda},\nabla u_k)dx\nonumber\\
& \le & - \int_{{\R^N}}V(x)\phi(|u_{\lambda}|)u_{\lambda}\varphi_kdx+\lambda \int_{{\R^N}} \widetilde{f}(x,u_{\lambda})\varphi_kdx\nonumber\\
& \le &\lambda D_4 \int_{{\R^N}} |u_{\lambda}|^{r}|u_k|^{l(\beta-1)}dx\nonumber\\
& =  &\lambda D_4 \int_{{\R^N}} |u_{\lambda}|^{r-l}|w_k|^{l}dx.\nonumber
\end{eqnarray}
The rest proof  is the same as Lemma 2.6 in \cite{Medeiros} with replacing $p$ with $l$ and $p^*$ with ${l}^*$.
\qed

 \vskip2mm
 \noindent
{\bf Proof of Theorem 5.1.}
 By Lemma 5.3 and Lemma 5.4, we have
\begin{eqnarray*}
          \|u_{\lambda}\|_{\infty}
  \le   C\lambda^{\frac{1}{l^*-r}}C_{**}^{\frac{l^*-l}{l^*-r}}
 \max\left\{\lambda^{-\frac{1}{k-l}\cdot{\frac{l^*-l}{l^*-r}}},\lambda^{-\frac{l}{m(k-l)}\cdot {\frac{l^*-l}{l^*-r}}}\right\}.
\end{eqnarray*}
 Notice that $k>l, l^*>l$, $l^*>r$ and
\begin{eqnarray*}
 k<K=\min\left\{l^{*},\frac{ml-l+l^*}{m}\right\}.
\end{eqnarray*}
 Then
there exists a large $\Lambda_0>0$ such that $ \|u_{\lambda}\|_{\infty}<\frac{\delta}{2}$  for all $\lambda>\Lambda_0$, which implies that
$\widetilde{F}(x,u_{\lambda})={F}(x,u_{\lambda})$ for all $x\in \R^N$. Hence, $u_{\lambda}$ is a nontrivial weak solution of system (\ref{741}) and Lemma 5.3 and Lemma 5.4 imply that  $\|u_\lambda\|_{1,\Phi}\to 0$ and   $\|u_\lambda\|_{\infty}\to 0$  as $\lambda\to\infty$, respectively.\qed
 \vskip2mm
 \noindent
 {\bf Remark 5.1.} Comparing Theorem 1.1 with Theorem 5.1, it is easy to see Theorem 5.1 for the scalar equation (\ref{741}) is better because $f$ satisfies growth conditions just for $|t|\le \delta$ with some positive constant $\delta$ rather than for $|t|\le 2$. The proof of Theorem 1.1 for the elliptic system (\ref{aa1}) present  more complex derivation. Especially, Moser iteration for system (\ref{aa1}) is  more difficulties than  that for the scalar equation (\ref{741}). Moreover, in Theorem 1.1, we assume that $F$ satisfies (F1)-(F3) for all $|(t,s)|\le 4$ which is a circular domain. However, if we assume that there exist two positive constants $\delta_1$ and $\delta_2$ such that $F$ satisfies (F1)-(F3) for all $|t|\le \delta_1$ and $|s|\le \delta_2$ which is a rectangular domain, the arguments will become more complex and it is unknown if Theorem 1.1 holds in such rectangular domain.  One can consider the proof of Lemma 3.1 and the examples of cut-off functions to see such complexity.   Finally, we would like to mention that if
 we let
  \begin{eqnarray}\label{aaa1}
 \widetilde{F}^+(u)=\begin{cases}
  \widetilde{F}(u) & \mbox{if } u\ge 0\\
  0 & \mbox{if } u < 0\\
   \end{cases}
  \  \mbox{,} \ \ \ \ \ \;\;\;\;\;
    \widetilde{F}^-(u)=\begin{cases}
  \widetilde{F}(u) & \mbox{if } u\le 0\\
  0 & \mbox{if } u > 0\\
   \end{cases}
 \end{eqnarray}
 and consider the functional
  $$
  \widetilde{J}_{\lambda}(u)=\int_{\R^N}\Phi(|\nabla u|)dx+\int_{\R^N}V(x)\Phi(|u|)dx-\lambda\int_{\R^N} \widetilde{F}^+(x,u)dx
  $$
  and
  $$
  \widetilde{J}_{\lambda}(u)=\int_{\R^N}\Phi(|\nabla u|)dx+\int_{\R^N}V(x)\Phi(|u|)dx-\lambda\int_{\R^N} \widetilde{F}^-(x,u)dx,
  $$
  respectively, then we can obtain equation (\ref{741-1}) has a positive solution and a negative solution. Thus Theorem 5.1 can be seen as the generalization of Theorem 1.2 in \cite{Medeiros} if we assume that $(\mathcal{V}1)$ holds instead of the periodicity of $V$.

\section*{Acknowledgments}
This work is supported by Yunnan Ten Thousand Talents Plan Young \& Elite Talents Project and Candidate Talents Training Fund of Yunnan Province (No: 2017HB016). This work was completed partially during the first author's visit to the Western University. The first author would like to thank the support of the China Scholarship Council for study abroad and the faculties in the Department of Applied Mathematics  of the Western University for their warm help.

 {}
 \end{document}